\documentclass[moor]{informs1}              

\usepackage{natbib}
 \NatBibNumeric
 \def\bibfont{\small}%
 \bibpunct[, ]{[}{]}{,}{n}{}{,}%

\usepackage[colorlinks=false,urlcolor=blue,citecolor=blue,linkcolor=blue,bookmarks=true,bookmarksopen=false,pdftitle=created by dvipdf,pdfcreator=NM,pdfauthor=Jain,pdfsubject=mor.sty:6.29.04]{hyperref}
\usepackage{amsmath,amstext,amsfonts,amssymb,mathrsfs,euscript}
\usepackage{graphicx,color,fullpage,setspace,subfigure,cases}

\usepackage{personal2}
\newcommand{\sn}{\sqrt{n}}

\newcommand{\tY} {\tilde{Y}}
\newcommand{\tT}{\tilde T}

\TITLE{On Transitory Queueing}
\ARTICLEAUTHORS{%
\AUTHOR{Harsha Honnappa}
\AFF{EE Department, University of Southern California\\ \EMAIL{\tt honnappa@usc.edu}}
\AUTHOR{Rahul Jain}
\AFF{EE \& ISE Departments, University of Southern California\\ \EMAIL{\tt rahul.jain@usc.edu}}
\AUTHOR{Amy R. Ward}
\AFF{Marshall School of Business, University of Southern California\\ \EMAIL{\tt amyward@usc.edu}}
}

\begin{document}

\ABSTRACT{
We introduce a framework and develop a theory of `transitory' queueing
models. These are models that are not only non-stationary and
time-varying but also have other features such as the queueing system
operates over finite time, or only a finite population arrives. Such
models are relevant in many real-world settings, from queues at
post-offices, DMV, concert halls and stadia to out-patient departments
at hospitals. We develop fluid and diffusion limits for a large class
of `transitory' queueing models. We then introduce three specific
models that fit within this framework, namely, the $\D_{(i)}/GI/1$
model, the conditioned $G/GI/1$ model, and an arrival model of
scheduled traffic with epoch uncertainty. We show that asymptotically
these models are distributionally equivalent, i.e., they have the same
fluid and diffusion limits. We note that our framework provides the
first ever way of analyzing the standard $G/GI/1$ model when we
condition on the number of arrivals. In obtaining these results, we
provide generalizations and extensions of the Glivenko-Cantelli and
Donsker's Theorems to triangular arrays. Our  analysis uses a
technique we call \textit{population acceleration}, which we discuss
in some detail. 
}


\MSCCLASS{Primary: 60K25, 90B15; secondary: 68M20, 90B22}
\ORMSCLASS{Primary: queues, diffusion models, limit theorems; secondary: nonstationary, transient results}
\HISTORY{This version: \today}

\maketitle

\section{Introduction}\label{sec:intro}

Erlang's study in 1909 \cite{Er1909}, of what came to be known as the $M/D/1$ queue, initiated the theoretical study of queueing phenomena. He argued that the number of calls arriving at a telephone exchange in a given time interval is Poisson distributed. 
Subsequent work in queueing theory generalized the Poisson traffic model to
renewal, or $GI$ for ``general independent'' inter-arrival time, traffic processes. This enlarged the class of models that could be studied, and as Kingman \cite{Ki2009} notes, while remaining largely tractable analytically. The early analyses of $GI/G/1$ queueing models
focused on stationary and ergodic queueing systems. Seminal work by Pollaczek, Kendall, Kingman (among others) developed a comprehensive theory of stationary and ergodic queueing systems  (see \cite{Ki2009} for a chronological perspective of this work.)

Ergodic analysis has been hugely important in revealing the ``large time'' (as $t \to \infty$) or steady state behavior of queueing systems. But often, we are interested in transient, or ``small time'', analysis of queueing systems. Typically,  this is quite messy even for simple models such as $M/M/1$ (see, for example, \cite{GrHa98}). One of the more celebrated results in this context is that of Iglehart and Whitt \cite{IgWh70} who showed that under the heavy traffic condition, the transient  queue length and waiting time are distributionally approximated by reflected Brownian motion processes. This is called the heavy traffic diffusion approximation \cite{Ha85,Gl1990}.  

In reality, queueing systems often exhibit non-stationarities in their operation, and heavy traffic conditions may not always pertain. For instance, in \cite{BrGa2005} the authors analyze  data at a call center, and show that the traffic is approximately like a non-homogeneous (time-varying) Poisson process \cite{KiWh2013b,KiWh2013a}. Much of the literature on non-stationary queues has focused on the Markovian case (or $M_t/M_t/\cdot$ queues). Newell \cite{Ne68a,Ne68b,Ne68c,Ne82} had the earliest, heuristic analysis of such systems. This was expanded upon  by \cite{Ke82,RoOr1979} via a \emph{pointwise stationary  approximation} (PSA) and formalized in  \cite{Ma81,MaWh1998,Ma85} via the \emph{uniform acceleration} (UA)  technique. The key idea behind UA is that by rescaling the  arrival and service rates of the system, one obtains a ``slowly  varying'' system whose mixing time is sufficiently small that it reaches stationarity quickly. However, many queueing systems are non-ergodic in nature.  This is either because the queue is unstable, or the system is varying rapidly. In such cases, UA/PSA can be a poor approximation.

Consider queues at post-offices, concert halls, stadiums, retail stores during black friday sales, or scheduled arrivals at a hospital out-patient department. In each instance, there is either a finite number of arriving users, or service is offered for a fixed period of time. These queueing systems are transient, non-stationary and non-ergodic in nature. We call such systems \emph{transitory queueing systems}. They are very difficult to analyze via classical queueing theoretic techniques since approximating the system state at any time $t$ by a notional stationary and ergodic process may not be tenable. In this paper, we introduce a general model of `transitory queues', specify a number of queueing models that fit within this framework, and introduce methods and techniques for their analysis. Exact analysis of transitory queueing models is almost impossible, and one must resort to approximations. Moreover, the UA technique  is hard to use as such models are non-Markovian, and the arrival and service rates need not vary slowly enough for even local ergodicity to hold.

In this paper, we use an alternative approach that we call the \emph{population acceleration} (PA) technique. In this technique, the queueing performance metrics are studied by increasing the number of users in a fixed time interval. 
The PA technique is similar to the UA technique in that the number of users arriving in a small interval goes to infinity but yet distinct, as the time axis is not scaled, and it doesn't depend on some form of ergodicity holding. We derive functional Strong Law of Large Numbers (fSLLN)/fluid limit and functional Central Limit Theorem (fCLT)/diffusion limit approximations to the queue length process as the population size increases. We consider a sequence of $s$-server transitory queues with FIFO service. Our only assumptions on the traffic model are Assumption 1: the fluid limit of the arrival process is a cumulative distribution function (i.e., right continuous and nondecreasing with finite total mass), and Assumption 2: the diffusion limit is a tied-down Gaussian process, possibly with continuous sample paths. Under these fairly weak assumptions, we show that the fluid limit is a Skorokhod reflection of the fluid limit of an appropriately defined netput process. It is, in general, time-varying, and switches between `overloaded', `underloaded' and `critically loaded' regimes. The diffusion limit of the queue length process is shown to be a \emph{directional derivative} of the one-dimensional Skorokhod reflection map of the netput fluid limit process in the direction of the diffusion netput process, which is a combination of a tied-down Gaussian process and an independent Brownian motion process. To study convergence, we introduce the space $\sD_{\lim}[0,\infty)$, the space of all functions that are right \textit{or} left continuous with right \textit{and} left limits  at every point, and right continuous at 0, which is larger than $\sD$ (the space of \emph{cadlag} functions). Moreover, convergence is obtained w.r.t. the $M_1$ topology on $\sD_{\lim}[0,\infty)$, and we provide a counterexample to show that convergence in the stronger $J_1$ is not possible for transitory queueing models in general.

We refer to queueing models that satisfy Assumptions 1 and 2 as \textit{transitory queueing models}. This is a fairly broad class. We show three traffic models that satisfy these assumptions. 
\begin{enumerate}
\item[(i)] \textit{The $\D_{(i)}$ Traffic model:} In the basic $\D_{(i)}$ model, we assume that the arrival times of users are sampled independently from an \textit{identical distribution}. Thus, the arrival times are ordered statistics and the inter-arrival times (denoted $\D_{(i)}:=T_{(i+1)}-T_{(i)}$) are differences of ordered statistics, hence the name the  $\D_{(i)}$ traffic model. This model is studied in detail in \cite{HoJaWa12} where we develop fluid and diffusion limits for the queue length and a transient Little's law.

Here, we consider the \textit{generalized} $\D_{(i)}$ traffic model wherein the arrival time of each user is independently sampled from \textit{non-identical distributions}. Thus, the unordered arrival times now form a triangular array. To study this model, we provide generalizations of the Glivenko-Cantelli and Donsker's Theorem for triangular arrays. We do so via a  generalization of Hahn's Central Limit Theorem \cite{Ha1978} to non-identically distributed processes. Using the notion of a \emph{random distribution function} introduced by Dubins and Freedman \cite{DuFr1967} we identify the fluid and diffusion limits for the generalized theorems, thus proving that the \textit{generalized} $\D_{(i)}$ traffic model satisfies Assumptions 1 and 2.

\item[(ii)] \textit{The Conditioned Renewal Process model:} The $\D_{(i)}$ traffic model seems very natural, particularly to those uninitiated in queueing theory. And yet, queueing theory has mostly focused on the renewal process traffic model. A key question is how are the two related, if at all?  Thus, we next introduce the conditioned renewal arrival process model. Herein, the arrivals happen according to a renewal process but we condition on there being $n$ arrivals by some time $T$. When the renewal process is Poisson, it is well known that the joint distribution of the arrival times is Uniform when conditioned (which is the $\D_{(i)}$ arrival model with the Uniform distribution). From this, it is easy  to conclude that the arrival processes are also equal in distribution. It is well known that a conditioned renewal process is not distributionally equivalent to a model with i.i.d. sampling from some distribution $F$. Herein, we show that, in fact, the conditioned renewal arrival process is asymptotically distributionally equivalent to a $\D_{(i)}$ arrival model with some distribution $F$ as $n \to \infty$, in the sense that both processes converge to the same weak limit, namely a Brownian bridge process.


\item[(iii)] \textit{Scheduled Arrivals with Epoch Uncertainty model:} Many queueing scenarios involve scheduled arrivals at appointment times, e.g., arrivals at a doctor's office. And yet, there is randomness in the actual arrival time, around the scheduled time. The earliest reference to such a model of traffic is in \cite{CoSm1961}, where it is referred to as ``a regular arrival process with unpunctuality.'' There is an increasing interest in such queueing models which are not amenable to analysis via known queueing theoretic methods. We model randomness in arrival times around the scheduled times  as being uniformly distributed over a small interval. We show that even though the sampling model is different from the  $\D_{(i)}$ model, this is also a transitory queueing model that satisfies Assumptions 1 and 2, and has the same weak limit as the $\D_{(i)}/GI/1$ queueing model. 

\end{enumerate}

This confluence of asymptotes for some disparate but natural models for `transitory' queueing phenomena is an interesting coincidence. In fact, one may potentially construct other models as well that satisfy Assumptions 1 and 2. In all such cases, the fluid and diffusion limits will be the same as for the $\D_{(i)}/GI/s$ queueing model. It is worth mentioning here that the $\D_{(i)}/GI/1$ model arose as an equilibrium model in \cite{HoJa11,JaJuSh11} where a finite population of users were considered to be strategically picking their arrival times. Thus, in some sense the $\D_{(i)}/GI/1$ queueing model can be considered canonical to the study of transitory queues just as the $M/GI/1$ and $G/G/1$ queueing models are to the study of stationary queueing systems.


\subsection*{Literature Review.} \label{sec:lit-rev}

There has been a long-standing interest in non-stationary and time-varying queueing models. One of the first pieces of work is that of Newell \cite{Ne68a,Ne68b,Ne68c} who characterized the various operating states of the non-homogeneous Poisson queue as the load factor $\rho(\cdot)$  varies with time. The motivation came from transportation networks, and Newell performed a heuristic analysis. Keller \cite{Ke82} provided more formal arguments and showed  that the transient distribution  at time $\t$ can be  approximated by the stationary measure associated with a notional  Markov chain that has arrival and service rates $\l(\t)$ and  $\m(\t)$, respectively. This type of analysis has come to be known as the  \emph{pointwise stationary approximation} (PSA) (see \cite{Wh1991}  as well). Massey \cite{Ma81,Ma85} and Massey and Whitt \cite{MaWh1998} made these arguments more rigorous and showed that  Keller's perturbation approach can be justified as a \textit{uniform  acceleration} (UA) asymptotic expansion of the transient  distribution. The notion of uniform acceleration comes from the fact  that the arrival and service rates are scaled at all time instants by the same parameter $\e$, and the expansions arise as $\e \to 0$. Later, Mandelbaum and Massey \cite{MaMa95} developed  fSLLN and fCLT results using Strong Approximation techniques for the $M_t/M_t/1$ queue, and identified the directional derivative reflection map as the right one to succinctly represent the queue length process diffusion limit in all regimes. This is based on the UA technique which relies on the assumption that the time scales on which the queue can change appreciably is of the order of $1/\e$ for some $\e > 0$. This technique has been extensively applied to non-stationary queueing systems with non-homogeneous Poisson input \cite{Ma2002}. However, it is not yet clear whether it is also useful for transitory queueing models of the kind we introduce in this paper.  

More closely related is the ``Binomial Traffic Model'' that Newell \cite{Ne82} introduced. This corresponds to the i.i.d. sampling $\D_{(i)}$ model. Through heuristic analysis, Newell identifies the limit processes in different regimes. However, these are point-wise and not functional limits, and a weak convergence result is missing. In \cite{GaLePe75}, the authors also identify several practical scenarios where an i.i.d. sampled $\D_{(i)}$ traffic model would be meaningful. Some analysis is also presented when the arrival process is approximated to be Markovian. The birth-death transient balance equations are solved numerically, and it is shown that a ``deterministic approximation'' (i.e., a first order fluid model) is good as the population size  increases. We, however, note the the main difficulty in analyzing the $\D_{(i)}/GI/1$ model is precisely the lack of Markovian structure. The  i.i.d. sampling $\D_{(i)}$ traffic model was also studied by  Louchard \cite{Lo94}. He provided an analysis analogous to Newell \cite{Ne68a,Ne68b,Ne68c} in the time varying Markovian queue case, and established the diffusion limits at continuity points in  certain regimes (over-saturated and near-saturation). However, these are not functional limits. Moreover, the whole difficulty in establishing diffusion limit for the $\D_{(i)}/GI/1$ model is precisely the fact there are discontinuities as the limit process switches regimes. 

Thus, we see that there has been long-standing interest in non-stationary, time-varying queueing models. In fact, there has even been an interest in modeling `transitory' queueing phenomena. In this paper, we provide a framework and a class of `transitory' queueing models. We show the connections between three interesting, and somewhat natural `transitory' queueing models. Moreover, we establish fluid and diffusion limits for the whole class of such models. We do this by using the \textit{population acceleration} technique as opposed to the uniform acceleration technique that has been usually used in studying non-stationary queues. In establishing these results, we have had to establish or generalize existing results on Glivenko-Cantelli and Donsker's Theorems for empirical processes to triangular arrays. Those mathematical results should be of independent interest, and we hope will be useful elsewhere as well.

The paper is organized as follows. Section \ref{sec:prelim}  introduces and defines a general transitory queueing model. Section \ref{sec:performance-analysis} develops fluid and diffusion limits for performance metrics of transitory queueing models that satisfy Assumptions 1 and 2. Section \ref{sec:traffic-models} introduces three transitory queueing models that satisfy Assumptions 1 and 2, and show that the fluid and diffusion limits for all three coincide. Conclusions and discussion about future work is provided in Section \ref{sec:conclusion}.

\section{The Transitory Queueing Model}\label{sec:prelim}

There is a finite population of $N$ customers that arrive to an infinite buffer for service.  The service opens at time 0; however, some customers arrive beforehand.  The earliest possible arrival time is $-T_0 \leq 0$.  The random vector $\mathbf T := (T_{1}, T_{2}, \cdots, T_{N}) \in [-T_0,\infty)^N$ represents the arrival times of the $N$ customers.  We assume all elements of $T$ are finite with probability 1.  The cumulative number of arrivals up to time $t$ is
\begin{equation} \label{eq:arrival-process}
{A}(t) =  \sum_{i=1}^n  \mathbf{1}_{\{T_{i} \leq t\}}.
\end{equation}
We call the ${A}(0)$ customers that arrive before the service opens {\em early birds}.  

The $s$ servers process the arrivals in a first-come-first-served manner.   The servers are non-idling and service is non-preemptive.  The $i$th customer to receive service from server $j \in \{1, \ldots, s\}$ has processing time $\nu_{j,i}$, which has CDF $G$ and support $[0,\infty)$.
The $s$ i.i.d. infinite sequences of processing times
$\{\nu_{j,i}, i \geq 1\}$ are mutually independent of each other and of the arrival epochs $\mathbf{T}$.  The service time mean is $1/\mu := \bbE\nu_{j,i}  <
\infty$, and the variance is $\sigma^2 := \mbox{Var}(\nu_{j,i}) < \infty$.  The number of potential service completions if server $j$ was busy in all of $[0,t]$ is given by the renewal counting process $S_j$, defined as
\begin{equation} \label{def:service-process}
S_j(t) := \sup \{m \geq 1 | V_j(m) \leq t\}, ~\forall t \in [0,\infty),
\end{equation}
where
\(
V_j(m) := \sum_{i=1}^m \nu_{j,i}.
\)



Now, let $Q$ represent the \emph{queue length}
process, including any customers in service and all buffered
customers. The sample paths of $Q$ are defined in terms of those of the arrival and service processes as
\begin{equation} \label{def:queue-length}
Q(t) := A(t) - \sum_{j=1}^s S_j(B_j(t)) \geq 0, \quad \forall t \in [-T_0,\infty),
\end{equation}
where $B_j(t)$ is the \emph{busy time} of server $j$, defined as the amount
of time server $j$ spent serving jobs in the interval $[0,t]$. Let
$B(t) := \sum_{j=1}^s B_j(t)$ be the total busy time of the queue.
 
When $s=1$, it follows that $B_1(t) := \int_0^t \mathbf{1} \{Q(s) >0 \} ds$, for all $t>0$; however, in general the characterization of each $B_j$ is complex, and depends on how arriving customers that find more than one server idle are routed.  We do not provide an explicit representation for $B_j$.  Instead we provide conditions that must be satisfied by restricting when servers can be idle.  Our analysis applies to any routing policy that satisfies those conditions.  For a concrete example, the reader may assume that when an arriving customer finds more than one server idle, that arrival is equally likely to be served by any of the idle servers.

The cumulative idle time of server $j \in \{1,2,\ldots,s\}$ in $[0,t]$ is
\[
I_j(t) := t \mathbf{1}_{\{t \geq 0\}} - B_j(t) \quad \forall t \in [-T_0,\infty).
\]
Note that it is natural to track only how much idle time each server has had since the service opened at time 0, despite the fact that customers may have been waiting before time 0, when the servers were ``off-duty''.  The total cumulative idle time of all the servers
\[
I(t) := \sum_{j=1}^s I_j(t)
\]
must satisfy
\(
I(0) = 0, \; I \mbox{ is non-decreasing, and } I(t) \mbox{ increases only if } Q(t) \leq s.
\)

Our {\bf objective} is to characterize the time-dependent queue-length distribution in our transitory queueing model  However, the queue-length process
is non-Markovian in general, which makes the analysis very difficult.  Moreover, even in the special case of a Markovian transitory queueing model, exact analysis does not result in closed-form expressions for the transient queue-length distribution, as the example below shows.
Our approach is to develop asymptotic approximations for the queue-length process as the population size $N$ becomes large.

\subsection*{The $\Delta_i$/M/1 Queue}
A special case to keep in mind is
when the joint distribution on $\mathbf T$ is product form, and the
marginals are identically distributed per a
given distribution function $F$ having density function $f$.  This is the model introduced in Newell~\cite{Ne82}.
Suppose also that the service times are exponential with rate $\m$.
Then, the joint variable
$(Q(t), \mathcal S(t))$, where $\mathcal S(t) \in 2^{\{1,\ldots,N\}}$ is the subset of
users who have arrived by time $t$, is a inhomogeneous continuous time Markov chain. By standard
`birth-death' chain arguments, observe that the transient state probability distribution
evolves per the following ordinary differential equation:
\[
\frac{d P(t, m,\mathcal{S})}{dt} =
\begin{cases}
\begin{split} &- \bigg( (N - | \mathcal{S}|) \frac{f(t)}{1-F(t)} \bigg ) P(t,m, \mathcal{S})\\ &+ |\mathcal{S}| \frac{f(t)}{1-F(t)} \sum_{i \in \mathcal{S}} P(t,m-1, \mathcal{S} \backslash \{i\}) \end{split} \qquad &\text{ if } t \leq 0,\,\, 0 < m \leq N \\ \\ \\
\begin{split} - \bigg( \mu &+ (N - | \mathcal{S}|) \frac{f(t)}{1-F(t)} \bigg ) P(t,m, \mathcal{S}) + \mu P(t,m+1,\mathcal{S}) \\ &+  |\mathcal{S}| \frac{f(t)}{1-F(t)} \sum_{i \in \mathcal{S}} P(t,m-1, \mathcal{S} \backslash \{i\}) \end{split} \qquad &\text{ if } t > 0,\,\, 0 < m \leq N
\end{cases}
\]
where $P(t,m, \mathcal S) = \bbP(Q(t) = m, \mathcal S(t) = \mathcal
S)$ and $t \in \bbR$.
These first-order differential equations are not easy to solve analytically.  Thus, even the simplest transitory queueing model is not amenable to exact analysis.

\subsection*{Notation}\label{sec:notation}
Unless noted otherwise, all intervals of time are subsets of
$[-T_0,\infty)$, for a given $-T_0 \leq 0$. Let $\sD_{\lim} :=
\sD_{\lim}[-T_0,\infty)$ be the space of functions $x : [-T_0,\infty)
\rightarrow \mathbb{R}$ that are right-continuous at $-T_0$, with
right and left limits and are either right or left continuous at every point $t > -T_0$. Note that
this differs from the usual definition of the space $\sD$ as the space
of functions that are right continuous with left limits (cadl\'{a}g
functions), and $\sD \subset \sD_{\lim}$. We denote almost sure convergence by
$\stackrel{a.s.}{\longrightarrow}$ and weak convergence by
$\Rightarrow$. The topology of convergence is indicated by the tuple
$(S,m)$, where $S$ is the metric space of interest and $m$ is the
metric that topologizes $S$. Thus, $X_n \stackrel{a.s.}{\longrightarrow}
X$ in $(\sD_{\lim},U)$ as $n \rightarrow \infty$ indicates that $X_n
\in \sD_{\lim}$ converges to $X \in \sD_{\lim}$ uniformly on compact
sets (u.o.c.) of $[-T_0,\infty)$ almost surely. Similarly, $X_n
\Rightarrow X$ in $(\sD_{\lim},U)$ as $n \rightarrow \infty$ indicates
that $X_n \in \sD_{\lim}$ converges weakly to $X \in \sD_{\lim}$
uniformly on compact sets of $[-T_0,\infty)$. $(\sD_{\lim},M_1)$
indicates that the topology of convergence is the $M_1$
topology.  We use $\circ$ to
denote the composition of functions or processes. The indicator
function is denoted by $\mathbf{1}_{ \{ \cdot \} }$ and the positive
part operator by $(\cdot)_+$. Finally, all random elements are defined
with respect to the canonical space $(\Omega, \sF, \bbP)$, unless noted
otherwise.\\

\section{Performance Analysis of Transitory Queueing Systems} \label{sec:performance-analysis}
In this section we present an analysis of the queue length performance
metric of the queueing model presented in Section
\ref{sec:prelim}, using the \emph{population acceleration}
technique. In Section \ref{sec:basic} we first establish the large
population asymptotics for the arrival and service processes for a generic
\emph{transitory} queueing model. Next, we establish fluid
approximations for the queue length by proving a fSLLN Theorem in
Section \ref{sec:fluid}. Finally, we establish a fCLT for the queue
length process and discuss its implications by considering a special
case in Section \ref{sec:diffusion} .

\subsection{Large Population Asymptotics of Primitives}\label{sec:basic}
We consider a sequence of systems indexed by $n \in \bbN$.  The customer population size in the $n$th system is $N_n$, and we assume $N_n \rightarrow \infty$ as $n \rightarrow \infty$.  Our convention is to superscript all processes and quantities associated with the $n$th system by $n$.

The arrival times in the $n$th system are $\mathbf T^n := (T_{1}^n, T_{2}^n, \cdots, T_{N_n}^n)$, and the cumulative arrival process $A^n$ is as defined in (\ref{eq:arrival-process}) with $T_i^n$ replacing $T_i$.  Our analysis requires that the empirical arrival distribution $A^n/N_n$ is well-behaved as $n$ becomes large.
\begin{assumption} \label{assume:arrival-process}
~There exists a probability distribution function $\overline{F}$ that has compact support such that the following holds.
\begin{enumerate}
\item[(a)] The arrival process satisfies a functional Strong Law of Large Numbers:
\[
\overline{A}^{n} := \frac{A^n}{N_n}  \stackrel{a.s.}{\to} \bar F \,\,
\text{in } (\sD_{\lim},J_1), \mbox{ as } n \rightarrow \infty.
\]
\item[(b)] The arrival process satisfies a functional Central Limit Theorem: 
\[
\hat A^n := \sqrt{N_n} \left( \frac{A^n}{N_n} - \bar F \right)
\Rightarrow \tilde{W} \,\,\text{in } (\sD_{\lim},J_1)  \mbox{ as } n \rightarrow \infty,
\]
\end{enumerate}
where $\tilde{W}$ is a zero mean Gaussian process with known
covariance function, that is tied-down to 0 at $-T_0$ and  $T_1 :=
\inf\{t: \bar F(t) = 1\}$. 
\end{assumption}
For example, when $N_n = n$ and $(T_1^n, \ldots,T_n^n)$ are
i.i.d. samples from a uniform distribution on $[0,1]$, the
Glivenko-Cantelli theorem guarantees that
Assumption~\ref{assume:arrival-process}(a) holds with $\overline{F}(t)
= t$ for $t \in [0,1]$, and, from Donsker,
Assumption~\ref{assume:arrival-process}(b) holds with $\tilde{W}$ a
standard Brownian Bridge (see \cite{BoSa96,KaSh91} for a formal
definition of a Brownian Bridge process and Theorem 13.1
in~\cite{Bi68} for a statement of Donsker's result).  However,
Assumption~\ref{assume:arrival-process} is also satisfied in much
greater generality, and we explore this in
Section~\ref{sec:traffic-models}.  In particular,
Assumption~\ref{assume:arrival-process} holds when arrival times are
sampled from different distributions, when $M^n$ is random, for a
conditioned renewal arrival model, and for a scheduled arrival
model. In all the models we have investigated, the limit $\tilde W$
is concentrated on $\sD \subset \sD_{\lim}$. With a small loss of
generality, we assume that $\bbP(\tilde W \in \sD) = 1$ for the
remainder of this paper.

The service times in the $n$th system are small; specifically, the service times of the $i$th arrival to server $j$ in the $n$th system is
\[
\nu_{j,i}^n := \frac{\nu_{j,i}}{N_n}, \; i=1,2,\ldots, \mbox{ for each } j \in \{1,\ldots,s\},
\]
so that $V_j^n(m) = \sum_{i=1}^m \nu_{j,i}/N_n$, and (\ref{def:service-process}) defines $S_j^n$ with $V_j^n$ replacing $V_j$.
 Furthermore, the fluid-scaled service process is
\[
\overline{S}^n(t) := \frac{1}{N_n} \sum_{j=1}^s S_j^n(t), \; t \geq 0,
\]
and the diffusion-scaled service process is
\[
\hat{S}^n(t) := \sqrt{N_n} \left( \overline{S}^n(t) - s \mu t \right), \; t \geq 0.
\]

Our analysis requires that the arrival and service processes, when
appropriately scaled, jointly satisfy a functional strong law of large
numbers and a functional central limit theorem. The multi-dimensional
result is shown to hold in the \emph{weak} $J_1$ topology $WJ_1$ on
$\sD_{\lim} \times \sD_{\lim}$; see Chapter 11 of \cite{Wh01} for details.
 
\begin{proposition} \label{proposition:basics}
~As $n \rightarrow \infty$, the fluid-scaled arrival and service processes jointly satisfy
\begin{equation} \label{lim:fluid-ASV}
(\bar{A}^n(t),\bar{S}^n(t) \mathbf{1}_{t \geq 0} )
\stackrel{a.s.}{\longrightarrow} (\bar F(t), s\mu t \mathbf{1}_{\{t
  \geq 0\}}) \text{ in } (\sD_{\lim} \times \sD_{\lim},WJ_1),
\end{equation}
and the diffusion-scaled arrival and service processes jointly satisfy
\begin{eqnarray}
\label{lim:diffusion-ASV}
(\hat{A}^n,\hat{S}^n) \Rightarrow \bigg( \tilde W, W \bigg) \text{ in
} (\sD_{\lim} \times \sD_{\lim},WJ_1),
\end{eqnarray}
where
\[
W(t) = \begin{cases}
\sigma \m^{3/2} \sum_{j=1}^s W_j(t) & ~~ t \geq 0 \\
0 & ~~ t < 0
\end{cases}
\]
is the sum of
independent standard Brownian
motion processes $W_j$, jointly independent of $\tilde W$ and $e: [0,\infty) \rightarrow [0,\infty)$ is the identity map.
\end{proposition}

 The proof of Proposition~\ref{proposition:basics} (see below) follows from Assumption~\ref{assume:arrival-process} and standard results on renewal processes, except for the subtlety that those results are usually proved in $\s D$ instead of $\s D_{\lim}$.  In particular, the following technical Lemma, used repeatedly throughout this paper, is useful to show (\ref{lim:diffusion-ASV}).  Its proof can be found in the appendix.

\begin{lemma} [Technical Lemma] \label{lem:d-lim}   Let $\mathscr{D}_{\lim}$ and $\mathscr{D}$ represent the
Borel $\s$-algebra generated by the $J_1$
topology on $\sD_{\lim}$ and $\sD$ (resp.)
\noindent (i) ~Let $x$ be a random element taking values in
the space $\sD$, where $\sD \subset
\sD_{\lim}$. Then, the measure induced by $x$ on $(\sD,\mathscr{D})$
can be extended to $(\sD_{\lim},\mathscr{D}_{\lim})$.

\noindent (ii) ~Let $\{x_n\}, ~n\geq 1$ be a collection of random
elements in $\sD$, such that $x_n \Rightarrow x ~\text{in}~(\sD,J_1)$
as $n \to \infty$. Then,
\(
x_n \Rightarrow x ~\text{in}~(\sD_{\lim},J_1)
\)
~as $n \to \infty$.
\end{lemma}

\Proof [Proposition \ref{proposition:basics}]
First note that by Assumption \ref{assume:arrival-process} and the assumed independence of the arrival and service processes, it is enough to show the convergence of $\overline{S}^n(t) \mathbf{1}_{t \geq 0}$ and $\hat{S}^n$.  These convergences hold in $(\sD, J_1)$ by the functional strong law and functional central limit theorems for renewal processes (see, for example,  Chapter 5 of \cite{ChYa01} and Theorem 7.3.2 and Corollary 7.3.1 in \cite{Wh01}), and the continuity of the addition operator with respect to the $J_1$ topology when the processes are continuous.  Finally, the convergence in $(\sD_{\lim}, J_1)$ in (\ref{lim:fluid-ASV}) is immediate since the measure induced by the limits concentrates degenerately on the fixed sample paths of the limit process in  $\sD \subset \sD_{\lim}$, and in (\ref{lim:diffusion-ASV}) is immediate from Lemma~\ref{lem:d-lim}.
\EndProof

The transitory queueing system model having customer population size $N_n$ (in the nth system) is defined as in Section 2, except that $A^n$ replaces A in (1) and $S^n$ replaces S in (2).  Then, the queue-length process $Q^n$ evolves as in (3) with $A^n$ and $S^n$ replacing A and S, and the busy time process $B^n$ replacing B.  The busy time process $B^n$ is defined through the idle time process $I^n$ accordingly.

\subsection{Fluid Approximations} \label{sec:fluid}
 We first derive the fluid limit for the queue-length process, and
 then the limit to the busy time process. Recall the definition of the
 queue length process in \eqref{def:queue-length}. The \emph{fluid-scaled} queue length process is
\begin{equation} \label{queue-length-fluid-scaled}
\bar Q^n(t) := \frac{Q^n(t)}{N_n} \,=\, \frac{1}{N_n} A^n(t) - \frac{1}{N_n} \sum_{j=1}^s S_j^n(B_j^n(t)),
\end{equation}
where $B_j^n(t)$ is server $j$'s fluid-scaled busy time process. 
Centering the right hand side of \eqref{queue-length-fluid-scaled} by
adding and subtracting the corresponding fluid-scaled processes,
and introducing the process $(s\m t ~ \forall t \geq 0)$ we obtain
\begin{equation*}
\bar Q^n(t) = \bigg( \frac{A^n(t)}{N_n} - \bar F(t) \bigg) - \sum_{j=1}^s \bigg(\frac{S_j^n(B_j^n(t))}{N_n} -
\mu B_j^n(t) \bigg) + \bigg( \bar F(t) - s\mu t \mathbf{1}_{\{t \geq 0\}} \bigg) + \sum_{j=1}^s \mu I_j^n(t),
\end{equation*}
where
\(
I_j^n(t)=t\mathbf{1}_{\{t \geq 0\}}-B_j^n(t)
 \)
~is the fluid-scaled idle time process. $\bar Q^n$ is equivalently written as
\(
\bar{Q}^n(t) \,=\, \bar{X}^n(t) + \bar Y^n(t), \quad \forall t \in
[-T_0,\infty),
\)
~where
\[
 \bar{X}^n(t) := \bigg( \frac{A^n(t)}{N_n} - \bar F(t) \bigg) - \sum_{j=1}^s \bigg(\frac{S_j^n(B_j^n(t))}{N_n} -
\mu B_j^n(t) \bigg) + \bigg( \bar F(t) - s\mu t \mathbf{1}_{\{t \geq 0\}} \bigg)
\]
and 
\[
\bar{Y}^n(t) := \sum_{j=1}^s \mu I_j^n(t) = \m I^n(t).
\]

In preparation for the main theorem in this Section, recall that the
one-dimensional Skorokhod reflection map is a (Lipschitz) continuous functional under
the uniform metric, $(\Phi, \Psi) : \sD_{\lim} \to \sD_{\lim} \times \sD_{\lim}$ defined as
\[
x \mapsto \Psi(x)(t) := \sup_{-T_0 \leq s \leq t} (-x(s))_+ ~~ \forall t \in \bbR
\]
and
\[
x \mapsto \Phi(x)(t) := x(t) + \Psi(x)(t), ~~ \forall x \in \sD_{\lim}
~\text{and}~~\forall t \in \bbR
\]
The Skorokhod reflection map satisfies the following properties. The
proof of claim (i) is part of Theorem 3.1 of \cite{MaRa10}, while (ii)
is a standard property of the Skorokhod reflection map and a proof can
be found in \cite{ChYa01,Wh01}.

\begin{proposition} \label{prop:srm-props}
~\noindent (i) $\Psi(\cdot)$ is continuous with respect to the uniform topology on $\sD_{\lim}$.\\
~\noindent (ii) $\Psi(x)(t)$ is non-decreasing in $t$. Further, for any other
pair of processes $(z,y) \in (\sD_{\lim}, \sD_{\lim})$ such that $z =
x + y \geq 0$, $y$ is non-decreasing, $y(0) = 0$ and $y$ increases
only if $z(t) \leq k$, for $k > 0$, the following relations hold:
$ \Phi(x-k)(t) \geq z(t) - k \geq \Phi(x)(t) - k$ and $\Psi(x - k)(t)  \geq y(t) \geq \Psi(x)(t)$.
\end{proposition}
Clearly, if $y$ is such that it does not increase when $z > 0$, then
the inequalities in $(ii)$ match each other. This is called the
\emph{dynamic complementarity} property of the  one-dimensional Skorokhod reflection
map. Therefore, $(ii)$ defines an approximate dynamic complementarity property.
\begin{theorem} [Fluid Limit]
\label{thm:queue-length-fluid}
~The pair $(\bar{Q}^n, \bar{Y}^n)$ jointly converges as $n \rightarrow \infty$,
\[
(\bar{Q}^n, \bar Y^n) \stackrel{a.s.}{\longrightarrow} (\Phi(\bar{X}),
\Psi(\bar{X}))~ \text{ in } (\sD_{\lim} \times \sD_{\lim},WJ_1),
\]
where $\bar{X}(t) = (\bar F(t) - s \mu t \mathbf{1}_{\{t \geq 0\}})$.
\end{theorem}
\Proof
First note that $\bar{Q}^n(t) \geq 0, ~ \forall t \in  [-T_0,\infty)$.  It is also true that $I^n(-T_0)  = 0$ and $d I^n(t) \geq 0, ~ \forall t \in [-T_0,\infty)$. 
From (ii) in Proposition \ref{prop:srm-props} it follows that $\bar Q^n(t) \geq \Phi(\bar X^n)(t)$ and $\bar Y^n(t) \geq \Psi(\bar X^n)(t)$.

By definition, $B_j^n(t) \leq t$ for all $j = 1, \ldots, N$ and from
\eqref{lim:fluid-ASV} in Proposition  \ref{proposition:basics}, it follows that
\(
\left | \sum_{j=1}^s \left( \frac{S_j^n \circ B_j^n}{n} - \mu B_j^n \right) \right |
\stackrel{a.s.}{\longrightarrow} 0~ \text{ in } (\sD_{\lim},J_1).
\)
~Therefore, by applying  \eqref{lim:fluid-ASV} in Proposition \ref{proposition:basics} to the arrival
process it follows that
\(
\bar{X}^n \stackrel{a.s.}{\longrightarrow} \bar{X}~ \text{ in } (\sD_{\lim},J_1).
\)
~As a consequence of the limit derived above and the continuity of the
reflection map from (i) of Proposition \ref{prop:srm-props} we have
\(
\liminf_{n \to \infty} (\bar{Q}^n, \bar Y^n) \geq \lim_{n \to \infty } (\Phi(\bar{X}^n), \Psi(\bar{X}^n)) = (\Phi(\bar{X}),
\Psi(\bar{X}))~\text{ in } (\sD_{\lim} \times \sD_{\lim},WJ_1) \,\, a.s.
\)
~Next, using the upper bound in (ii) of Proposition \ref{prop:srm-props} we have the relation
\(
(\bar Q^n, \bar Y^n) \leq (\Phi(\bar X^n - \frac{s}{N_n}) + \frac{s}{N_n}, \Psi(\bar X^n - \frac{s}{N_n})).
\)
~As $s$ is fixed, it is obvious from the continuity of the
reflection map that
\(
\limsup_{n \to \infty} (\bar Q^n, \bar Y^n) \leq (\Phi(\bar{X}),
\Psi(\bar{X})) ~~ a.s.~\text{ in } (\sD_{\lim} \times \sD_{\lim},WJ_1)
\)
~This concludes the proof.
\EndProof

\noindent \textbf{Remarks} 1. Theorem \ref{thm:queue-length-fluid} shows that the fluid limit of the queue length process is
\(
\bar{Q}(t) \,=\, (\bar F(t) - s \mu t \mathbf{1}_{\{ t \geq 0\}})
+ \sup_{-T_0 \leq p\leq t} (-(\bar F(p) - s \mu p \mathbf{1}_{\{p \geq
  0\}}))_{+}, ~ \forall t \in [-T_0,\infty).
\)
~$\bar{Q}$ can be interpreted as the sum of the fluid netput process
and the potential amount of fluid lost from the system. Suppose that service
started with some workload in the system at time $0$ and that $(\bar
F(t) - s\mu t \mathbf{1}_{\{t \geq 0\}}) < 0$ for $t > 0$, so that the fluid
service process has ``caught up'' and exceeded the cumulative amount
of fluid arrived in the system by time $t$ (for simplicity assume $t >
0$). Let $\bar f$ represent the density function associated with the
distribution function $\bar F$ (if $\bar F$ has a
discontinuity at some point $t$, then $f(t) := \frac{f(t-) +
  f(t+)}{2}$). Suppose $\bar f(t) - s\mu < 0$, implying that the netput process is decreasing at
$t$. In this case, $\sup_{-T_0 \leq p \leq t} (-(\bar F(p) - s\mu p
\mathbf{1}_{\{p \geq 0\}}))_{+} = -(\bar F(t) - s\mu t)$. This is the amount of extra fluid that could have been served, but is now lost. \vspace{3 pt}

\begin{figure}[t]
\centering
\includegraphics[scale=0.5]{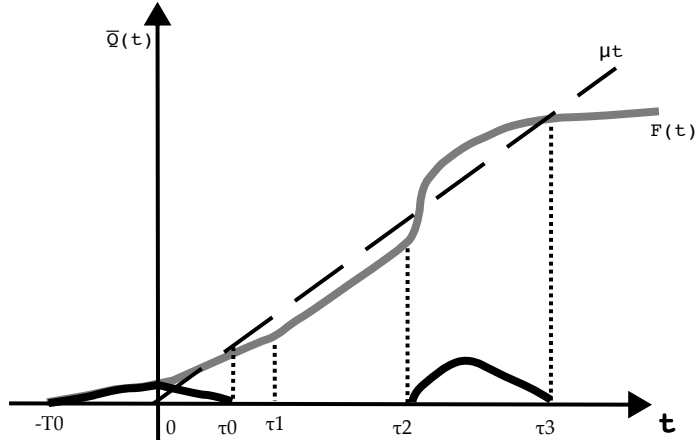}
\caption{An example of a fluid transitory queue that will undergo multiple
  ``regime changes''. The fluid queue length process is positive on $[-T_0,\tau_1)$ and $[\tau_2,\tau_3)$, and $0$ on
$[\tau_1,\tau_2)$ and $[\tau,\infty)$. Here $s = 1$.}
\label{fig:fluid}
\end{figure}
Figure \ref{fig:fluid} depicts an example queue length
process in the fluid limit for the special case in Section \ref{sec:special-case}, and its dependence on $F$ and
$\m$; here, $T = \t_3$ in the figure. In particular, notice that the process switches between being
positive and zero, during the time the queue operates. 
In particular, observe that these `regimes' correspond to when the
queue is overloaded, underloaded and critically loaded. It is
important to note that in any transitory queueing model, the (fluid)
limit system can experience these changes unlike a $G/G/1$
queue. Formally, these regimes can be codified by defining a `load
factor' $rho$ in terms of the fluid limit system as follows:

\begin{equation} \label{rho}
\rho(t) := \begin{cases}
\infty, \quad & \forall t \in [-T_0,0]\\
\sup_{0 \leq r \leq t}\frac{\bar F(t) - \bar F(r)}{\mu (t - r)}, \quad &\forall t \in[0,\tilde{T}]\\
0, \quad & \forall t > \tilde{T},
\end{cases}
\end{equation}
where $\tilde{T} := \inf \{t > 0 | \bar F(t) = 1 \text{ and } \bar{Q}(t) =
0\}$. Note that we define the traffic intensity to be $\infty$ in the
interval $[-T_0,0]$ as there is no service, but there can be fluid
arrivals. Based on this, we can now define the regimes of the 
transitory queueing model.

\begin{definition}[Operating regimes.]
~The transitory queue is
\begin{enumerate}
\item[\emph{(i)}] overloaded if $\rho(t) > 1$.
\item[\emph{(ii)}] critically loaded if $\rho(t) = 1$.
\item[\emph{(iii)}] underloaded if $\rho(t) < 1$.
\end{enumerate}
\end{definition}

This, in Figure \ref{fig:fluid} the queue is overloaded between
$[-T_0,\t_1]$ and $(\t_2,\t_3]$ and critically loaded between
$(\t_1,\t_2]$. It is possible to define finer states of the sytem, but we omit these
as they are not important to the central thesis of this
paper. However, the interested reader is encouraged to study the
analysis in \cite{HoJaWa12} where the sample paths of the diffusion
and fluid limits of the $\D_{(i)}/GI/1$ queue (a specific transitory
queueing model) are studied comprehensively. 

It is interesting to observe that the busy time of the queue, $B^n$,
does not converge to the identity process in contrast to the limit for
the GI/GI/1 queue in the heavy-traffic approximation setting. The
following corollary characterizes the busy time fluid limit.

\begin{corollary}
\label{cor:busy-time-fluid}
The fluid scaled busy time process $B^n := \sum_{j=1}^s B^n_j$
satisfies a fSLLN as $n \rightarrow \infty$:
\begin{equation}
\label{busy-time-fluid}
B^n \stackrel{a.s.}{\longrightarrow} \bar{B}  \text{ in } (\sD_{\lim},J_1)
\end{equation}
where $\bar{B}(t) \,:=\, s t \mathbf{1}_{\{t \geq 0\}} -
\frac{1}{s \mu} \Psi(\bar{X}(t))$, $\forall t \in [-T_0,\infty)$.
\end{corollary}
\Proof
By definition, we have $B^n(t) = s t \mathbf{1}_{\{t \geq 0\}} -
I^n(t) = s t \mathbf{1}_{\{t \geq 0\}} - \frac{\bar
  Y^n(t)}{\m}$. Theorem \ref{thm:queue-length-fluid} now implies the limit.
\EndProof

Note that $\bar{B}(t) = 0$ for all $t \leq 0$, as $\Psi(\bar{X})(t) =
0$ on that interval. It is important to keep in mind that $\bar B$ is
the \emph{total} busy time of the entire queueing system. For each
server, on the other hand, we can prove the following existence result.

\begin{corollary} \label{cor:server-busy-time-fluid}
~For every $j = 1,\ldots, s$, there exists a function $\bar B_j \in
\sC$ such that
\(
B^n_j \stackrel{a.s.}{\to} \bar B_j
\)
~ in $(\sD,J_1)$.
\end{corollary}
\Proof
Without loss of generality, let $t,~s~ \in [0,1]$. By definition,
${B^n_j}$ is a uniformly bounded sequence of functions (on the given
compacta) and $|B^n_j(t) - B^n_j(s)| \leq |t-s|$, for any such pair
$t,~s$. Thus, $B^n_j$ is uniformly Lipschitz, implying equicontinuity. Then, by the
Arzela-Ascoli Theorem the sequence $\{B^n_j\}$ is sequentially
compact, so that a limit exists.
\EndProof

 An exact characterization of $\bar B_j$ depends on the routing policy. However, this is not required for the rest of our analysis.

\subsection{Diffusion Approximations}\label{sec:diffusion}

Next, we derive the diffusion limit for the queue-length process in
Section~\ref{sec:queue-length-diffusion}. In Section
\ref{sec:special-case}, we specialize this result to a specific
instance of a transitory queueing system. Next, we show by a
counterexample by convergence in the $J_1$ topology is not possible in
general, in Section \ref{sec:M1J1}. Finally, we our main result to develop the diffusion limit for the busy time process. 

\subsubsection{Queue Length Process} \label{sec:queue-length-diffusion}
Define the \textit{diffusion-scaled} queue length process as
\begin{equation}
\frac{Q^n(t)}{\sqrt{N_n}} \,:=\, \frac{A^n(t)}{\sqrt{N_n}} -
\sum_{j=1}^s \frac{S_j^n(B_j^n(t))}{\sqrt{N_n}}, ~~\forall t \in [-T_0,\infty)
\end{equation}
Rewriting this expression by introducing the term $\sqrt{N_n} s\mu t
\mathbf{1}_{\{t \geq 0\}}$ and centering the terms on the right hand side
\begin{equation*}
\begin{split} \frac{Q^n(t)}{\sqrt{N_n}} = \bigg( \frac{A^n(t)}{\sqrt{N_n}} - \sqrt{N_n}
 \bar F(t) \bigg) - & \sum_{j=1}^s\bigg(\frac{S_j^n(B_j^n(t))}{\sqrt{N_n}} - \sqrt{N_n} \mu B_j^n(t) \bigg)\\
  &+ \sqrt{N_n} (\bar F(t) - s \m t \mathbf{1}_{\{t \geq 0\}}) + \sqrt{N_n} \sum_{j=1}^s \mu (t \mathbf{1}_{\{ t \geq 0\}} - B_j^n(t)).
  \end{split}
 \end{equation*}
Using the definition of the idle time process
\(
\sqrt{N_n} I_j^n(t) \,=\, \sqrt{N_n} (t\mathbf{1}_{\{ t \geq 0\}} - B_j^n(t)),
\)
~we can express $Q^n/\sqrt{N_n}$ as
\begin{equation} \label{def:queue-length-diffusion-scaled}
\frac{Q^n}{\sqrt{N_n}} = \hat{X}^n + \sqrt{N_n} \bar{X} + \hat{Y}^n
\end{equation}
where
\begin{eqnarray}
\label{X-hat-n}
\hat{X}^n(t) \,&:=&\, \bigg ( \frac{A^n(t)}{\sqrt{N_n}} - \sqrt{N_n}
\bar F(t) \bigg ) - \sum_{j=1}^s \bigg (\frac{S_j^n(B_j^n(t))}{\sqrt{N_n}} - \sqrt{N_n} \mu B_j^n(t)  \bigg )\\
\nonumber
&=&\, \hat{A}^n(t) - \sum_{j=1}^s \hat{S}_j^n(B_j^n(t)), \quad \forall t \in [-T_0,\infty),
\end{eqnarray}
and
\begin{equation}
\label{Y-hat-n}
\hat Y^n := \sqrt{N_n} \sum_{j=1}^s \mu I_j^n.
\end{equation}

Recall from Theorem \ref{thm:queue-length-fluid} that
$\bar{X}(t) =  (\bar F(t) -\m t \mathbf{1}_{t \geq 0}) $ is the fluid
netput process. We can think of $\hat{X}^n$ as a diffusion refinement
of the netput process. Lemma \ref{lem:X-hat} in the Appendix proves
that $\hat{X}^n$ converges weakly to a Gaussian process $\hat X$ as a direct consequence of
\eqref{lim:diffusion-ASV} in Proposition  \ref{proposition:basics}.

In the rest of this section, we will use Skorokhod's almost sure
representation theorem \cite{Sk56,Wh01b} and replace the random
processes above that converge in distribution by those defined on a
common probability space that have the same distribution as the
original processes and converge almost surely. The requirements for
the almost sure representation are mild; it is sufficient that the
underlying topological space is Polish (a separable and complete
metric space). We note without proof that the space $\sD_{\lim}$, as
defined in this paper, is Polish when endowed with the $J_1$
topology. This conclusion follows from Chapter 12.8 of \cite{Wh01} and the fact that
the proof there extends easily to $\sD_{\lim}$. The
authors in \cite{MaMa95} also point out that \cite{Po76} has a more
general proof of this fact. We conclude that we can replace the weak
convergence in \eqref{lim:diffusion-ASV} by
\begin{eqnarray*}
(\hat{A}^n,\hat{S}^n)
\stackrel{a.s.}{\longrightarrow} \bigg( \tilde W,W\bigg)~ \text{ in }
(\sD_{\lim} \times \sD_{\lim},WJ_1),
\end{eqnarray*}
where abusing notation we denote the new limit random processes by
the same letters as the old ones. This implies that in Lemma \ref{lem:X-hat}
\(
\hat{X}^n \stackrel{a.s.}{\longrightarrow} \hat{X}~ \text{ in }
(\sD_{\lim},J_1),
\)
~as $n \to \infty$.

Our goal is to establish diffusion limits for the centered queue length process
\begin{equation}
\label{queue-length-diffusion-centered}
\hat{Q}^n(t) \,:=\, \sqrt{N_n} \bigg (\frac{Q^n(t)}{N_n} - \bar{Q}(t)
\bigg ),
\end{equation}
and the process
\[
\tilde Y^n(t) := \hat Y^n(t) - \sqrt{N_n} \Psi(\bar X)(t).
\]
~We achieve this by using part (ii) of Proposition
\ref{prop:srm-props} to bound $(Q^n(t)/\sqrt{N_n}, \hat Y^n)$ in terms
of $\hat{X}^n$ and $\bar{X}$, and then establish the limit as $n \to
\infty$. The limit for each is proved in the weaker $M_1$ topology, as opposed
to the more common $U$ or $J_1$ topologies as convergence to the
\textit{directional derivative reflection map} (Lemma
\ref{lem:mama95} presented in the Appendix; see also \cite{HoJaWa12}
where a version of this theorem is proved in a special case) in general holds in $(\sD_{\lim}, M_1)$. In
fact, in Proposition \ref{thm:j1-convergence} below, a counterexample
is provided that shows that the limit result is not achievable in the
stronger $J_1$ topology, in general.

Recall that $(\Phi, \Psi)$ is the Skorokhod reflection map. The
directional derivative of the Skorokhod reflection map is defined below.

\begin{definition}[Directional Derivative Reflection Map] \label{def:dirdev}
~Let $x \in \sD$ and $y \in \sD$. For fixed $t \in [0,\infty)$
\begin{equation} \label{eq:dirderiv}
  \sup_{s \in \nabla_t^{x,L}}(-y(s-)) \vee \sup_{s \in \nabla_t^{x,R}}(-y(s)) := \lim_{a \rightarrow \infty} \Psi(a x + y)(t) - a \Psi(x)(t),
\end{equation}
is the directional derivative of $\Psi$ and
\(
\nabla_{t}^{x,L} := \{s \leq t | x(s-) = -\Psi(x)(t)\},
\)
~is a correspondence of points up to time $t$ where the \emph{left
  limits} of $x$ achieve an infinimum and
\(
\nabla_{t}^{x,R} := \{s \leq t | x(s+) = -\Psi(x)(t)\},
\)
~is a correspondence of points up to time $t$ where the \emph{right
  limits} of $x$ achieve an infinimum.
\end{definition}
Theorem 9.3.1 of \cite{Wh01b} proves the (pointwise) existence of the
limit. In establishing our main result, we use Theorem 9.5.1 of
\cite{Wh01b} and the Lipschitz continuity of the reflection map (in 
one-dimension) to prove the queue length diffusion limit in the $M_1$
topology. Of course, $M_1$ convergence is stronger than pointwise, in general.

\begin{theorem} [Diffusion Limit]
\label{thm:queue-length-diffusion}
\noindent (i) The diffusion scaled process $\tilde Y^n$ converges to a directional
derivative of the Skorokhod reflection regulator map:
\begin{equation} \label{directional-derivative-diff-limit}
\tilde Y^n \Rightarrow \tY ~\text{in}~ (\sD_{\lim},M_1)
\end{equation}
as $n \to \infty$, where $\tY(t) = \sup_{s\in
  \nabla_{t}^{\bar{X}, L}} (-\hat{X}(s-)) \vee \sup_{s\in
  \nabla_{t}^{\bar{X}, R}} (-\hat{X}(s))$ $\forall t \in
[-T_0,\infty)$ with $\nabla_t^{\bar X, \cdot}$ as in
Definition \ref{def:dirdev}.

\noindent (ii) The diffusion scaled queue length process $\hat{Q}^n$
converges to a reflected process, as $n \to \infty$:
\begin{equation} \label{queue-length-diffusion-limit}
\hat{Q}^n \Rightarrow \hat{X} + \tilde{Y}~ \text{ in } (\sD_{\lim},M_1),
\end{equation}
where $\hat{X}(t) = \tilde W(t) - \sum_{j=1}^s \sigma_j \mu_j^{3/2} W_j(\bar{B}_j(t))$.
\end{theorem}
\Proof
First, using \eqref{def:queue-length-diffusion-scaled} and the lower
bound in $(ii)$ of Proposition \ref{prop:srm-props} we have
\begin{eqnarray} \label{skorokhod-diffusion}
\bigg (\frac{Q^n}{\sqrt{N_n}}, \hat Y^n \bigg ) \,\geq\, \bigg ( \Phi(\hat{X}^n + \sqrt{N_n} \bar{X}), \Psi(\hat{X}^n + \sqrt{N_n}
\bar{X}) \bigg ).
\end{eqnarray}
This implies that
\begin{equation} \label{skorokhod-diffusion-centered}
\hat{Q}^n = \frac{Q^n}{\sqrt{N_n}} -  \sqrt{N_n} \bar{Q} \geq \Phi(\hat{X}^n + \sqrt{N_n} \bar{X}) - \sqrt{N_n} \bar{Q}.
\end{equation}

Recall from Theorem \ref{thm:queue-length-fluid} that $\bar{Q} =
\bar{X} + \Psi(\bar{X})$. Substituting this expression into
\eqref{skorokhod-diffusion-centered}, and using the fact that $\Phi(x)
= x + \Psi(x)$ for $x \in \sD_{\lim}$, we have
\begin{eqnarray}
\nonumber
\hat{Q}^n & \geq & \hat{X}^n + \sqrt{N_n} \bar{X} + \Psi(\hat{X}^n + \sqrt{N_n}
\bar{X}) - \sqrt{N_n} (\bar{X} + \Psi(\bar{X})),\\
\label{centered-queue-length}
 &=& \hat{X}^n + \Psi(\hat{X}^n + \sqrt{N_n} \bar{X}) - \sqrt{N_n}
 \Psi(\bar{X}).
\end{eqnarray}
Next, utilizing the expression for $\hat Y^n$ in
\eqref{skorokhod-diffusion}, and letting
\(
\tilde{Y}^n := \hat Y^n - \sqrt{N_n} \Psi(\bar{X}),
\)
~we have
\begin{equation}
\label{centered-idle-time}
\tilde{Y}^n \geq \Psi(\hat{X}^n + \sqrt{N_n} \bar{X}) - \sqrt{N_n} \Psi(\bar{X}).
\end{equation}
Therefore,
\begin{equation} \label{eq:diff-scale-lower}
\left( \hat{Q}^n, \tilde Y^n\right) \geq \left(\hat{X}^n + \tY^n, \Psi(\hat{X}^n + \sqrt{N_n} \bar{X}) - \sqrt{N_n} \Psi(\bar{X})\right).
\end{equation}

Next, using the upper bound in $(ii)$ of Proposition
\ref{prop:srm-props} we have
\[
\left( \frac{Q^n}{\sqrt{N_n}}, \hat Y^n \right)
\leq \left( \Phi \left(\hat{X}^n + \sqrt{N_n} \bar{X} -
  \frac{s}{\sqrt{N_n}} \right) + \frac{s}{\sqrt N_n}, \Psi\left(\hat{X}^n +
  \sqrt{N_n} \bar{X} - \frac{s}{\sqrt{N_n}} \right)\right).
\]
Now, using the centering arguments used in the lower bound we have
\begin{equation} \label{eq:diff-scale-upper}
\left( \hat Q^n, \tilde Y^n\right) \leq \left( \hat X^n + \Psi\left(\hat X^n + \sqrt{N_n} \bar X - \frac{s}{\sqrt{N_n}}\right) - \sqrt{N_n} \Psi(\bar X), \Psi\left(\hat X^n + \sqrt{N_n} \bar X - \frac{s}{\sqrt{N_n}}\right) - \sqrt{N_n} \Psi(\bar X) \right).
\end{equation}

The limit process follows by use of the directional derivative
reflection mapping lemma, Lemma \ref{lem:mama95} in the
Appendix. Using the fact that $\hat X^n \stackrel{a.s.}{\to} \hat X$
in $(\sD_{\lim},J_1)$, together with the lemma, it follows that
\(
\liminf_{n \to \infty} \tY_n \geq \tY := \sup_{s \in
  \nabla_{\cdot}^{\bar{X},L}} (-\hat{X}(s)) \vee \sup_{s \in
  \nabla_{\cdot}^{\bar{X},R}} (-\hat{X}(s)).
\)
~Consequently, from \eqref{eq:diff-scale-lower} we have
\(
\liminf_{n \to \infty} \hat Q^n \geq \hat X + \tY ~ \text{  in } (\sD_{\lim},M_1) \,\, a.s.
\)
~as $n \to \infty$.

Similarly, from \eqref{eq:diff-scale-upper}, and using Lemma
\ref{lem:X-hat} and Lemma \ref{lem:mama95} again, we have
\(
\limsup_{n \to \infty} \tilde Y^n \leq \tilde Y ~\text{in}~
(\sD_{\lim},M_1) ~~ a.s.
\)
~as $n \to \infty$. Then, using the upper bound on $\hat Q^n$ in
\eqref{eq:diff-scale-upper}, we have
\(
\limsup_{n \to \infty} \hat Q^n \leq \hat X + \tilde Y ~\text{in}~(\sD_{\lim},M_1) \,\, a.s.
\)
as $n \to \infty$. Combined with the lower bound above, this proves
convergence of the sample paths almost surely. Finally, the weak convergence in
 the statement of the theorem follows by the fact that the pre-limit processes
 are equal in distribution to our original processes.
\EndProof

\noindent \textbf{Remarks.}
1. Observe that the diffusion limit to the queue length process is a
function of a Gaussian bridge process and a Brownian motion
process. This is significantly different from the usual limits
obtained in a heavy-traffic or large population approximation to a
single server queue. For instance, in the $GI/G/1$ queue, one would
expect a reflected Brownian motion in the heavy-traffic setting.  In
\cite{MaMa95} it was shown that the diffusion limit process to the
queue length process of a $M_t/M_t/1$ queue is a time changed Brownian
motion $W(\int \l (s) ds + \int \m (s) ds)$, where $\l(s)$ and $\m(s)$ (resp.)
are the time inhomogeneous mean arrival rate and mean service rate
(resp.), reflected through the directional derivative reflection map
used in Lemma \ref{lem:mama95}. There are very few examples of
heavy-traffic limits involving a diffusion that is a function of a
Gaussian bridge and a Brownian  motion process; see example 3 of \cite{IgWh70}. However, there have been some results in other queueing models where a Brownian bridge arises in the limit. In \cite{PuRe10}, for instance, a Brownian bridge limit arises in the study of a many-server queue in the Halfin-Whitt regime. \vspace{3 pt}


\subsubsection{A Special Case} \label{sec:special-case}
To illustrate the difference between the population acceleration diffusion limit with the
RBM observed for the $GI/G/1$ queue, we present a corollary
to Theorem \ref{thm:queue-length-diffusion} when $\tilde W$ is a
Brownian Bridge process. A Brownian Bridge limit arises, for instance,
when the arrival times $T_{n,i}$ are sampled in an i.i.d. manner from
some distribution function $F$. We assume that $F$ has compact support
$[-T_0,T]$, where $T_0 > 0$ to allow for early bird arrivals, and that
the queue has a single server. This queue was comprehensively studied
in \cite{HoJaWa12}, where we call this model a $\D_{(i)}/GI/1$
queue. Notice that in this case, $\bar F = F$.

\begin{proposition} \label{prop:empirical-process-limit}
~If $T_{n,i}$ are i.i.d. samples from distribution function $F$, then
\[
\bar A^n \stackrel{a.s.}{\to} F ~~\text{in}~~ (\sD_{\lim},J_1)
\]
and
\[
\hat A^n \Rightarrow W^0 \circ F ~~\text{in}~~ (\sD_{\lim},J_1)
\]
as $n \to \infty$.
\end{proposition}
\Proof
We first fix $F$ to the uniform distribution function on $[0,1]$. The
fluid limit follows by the standard Glivenko-Cantelli; see Theorem
2.4.7 of \cite{Du10}. Theorem 16.4 of \cite{Bi68} proves that $\hat A^n \Rightarrow W^0~~\text{in}~~ (\sD, J_1)$. To prove that convergence holds in
$(\sD_{\lim},J_1)$, we must first show that the Brownian Bridge is well
defined on the larger space. But, this is a direct consequence of part
(i) of Lemma \ref{lem:d-lim}. Next, by part (ii) of Lemma
\ref{lem:d-lim}, it follows that $\hat A^n \Rightarrow
W^0~~\text{in}~~ (\sD_{\lim}, J_1)$. Finally, let $F$ be any arbitrary
cumulative distribution function. Since $W^0$ concentrates on
$\mathcal C \subset \sD$, Corollary 1 to Theorem 5.1 of \cite{Bi68} implies the
final result.
\EndProof

This proposition allows us to state a diffusion limit for the queue length
process of the $\D_{(i)}/GI/1$ queue, as a corollary of Theorem \ref{thm:queue-length-diffusion}.

\begin{corollary} \label{cor:brownian-bridge}
~Let $\tilde W = W^0 \circ F$ be a time changed Brownian Bridge
process. Then,
\begin{equation} \label{X-hat-integral}
\hat{X}(t) \stackrel{d}{=} \int_{-T_0}^t \sqrt{g^{'}(s)} d \check{W}_s, \quad \forall t \in [-T_0, \infty)
\end{equation}
where
\(
g(t) = F(t)(1-F(t)) + \sigma^2 \mu^3 \bar{B}(t)
\)
~and $\check{W}$ is a Brownian motion process on the same underlying
sample space $(\Omega, \sF, \bbP)$. Further, the queue length
diffusion limit process is
\begin{equation} \label{Q-Limit}
\hat Q(t) = \hat X(t) + \sup_{s \in \nabla_t^{\bar X}} (- \hat X(s))
~~ \forall t \in [-T_0,T],
\end{equation}
where $\nabla_t^{\bar X} := \{0 \leq s \leq t | \bar X(s) = -
\Psi(\bar X)(t)\}$ and $\bar X := F(t) - \m t$ is absolutely continuous.
\end{corollary}
\Proof
By Lemma \ref{lem:X-hat} it follows that $\hat X = W^0 \circ F - W
\circ \bar B$. By a classical time change (see, for
example, \cite{KaSh91}) $W^0 \circ F$ is equal in distribution to a
time changed Brownian motion, and $\hat{X}$ is equal in distribution
to the stochastic integral \eqref{X-hat-integral}. The diffusion function
$g(t)$ can be easily verified. The expression for $\hat Q$ now follows
by substitution in \eqref{queue-length-diffusion-limit}. Note that the
right and left correspondences $\nabla_t^{\bar X, R}$, $\nabla_t^{\bar
  X L}$  coincide since $\bar X$ is absolutely continuous.
\EndProof

\noindent \textbf{Remarks.} 1. Note that the diffusion $\hat X$ is
also tied down at some point in time, depending on whether $\t^* :=
\frac{1}{\m}$ is $> T$ or $\leq T$. In the former case, $\hat X(t) =
-\s \m^{3/2} W(\bar B(\t^*))$ for all $t \geq \t$, and in the latter
case $\hat X(t) = -\s \m^{3/2} W(\bar B(T))$ for $t \geq T$. Thus, the
process $\hat{X}$ can also be interpreted as  a time-changed Brownian
motion on the interval $[-T_0,T]$, tied down at a point determined by
the fluid busy time, with random tie-down level. \vspace{3 pt}

\begin{figure}[t]
\centering
\includegraphics[scale=0.6]{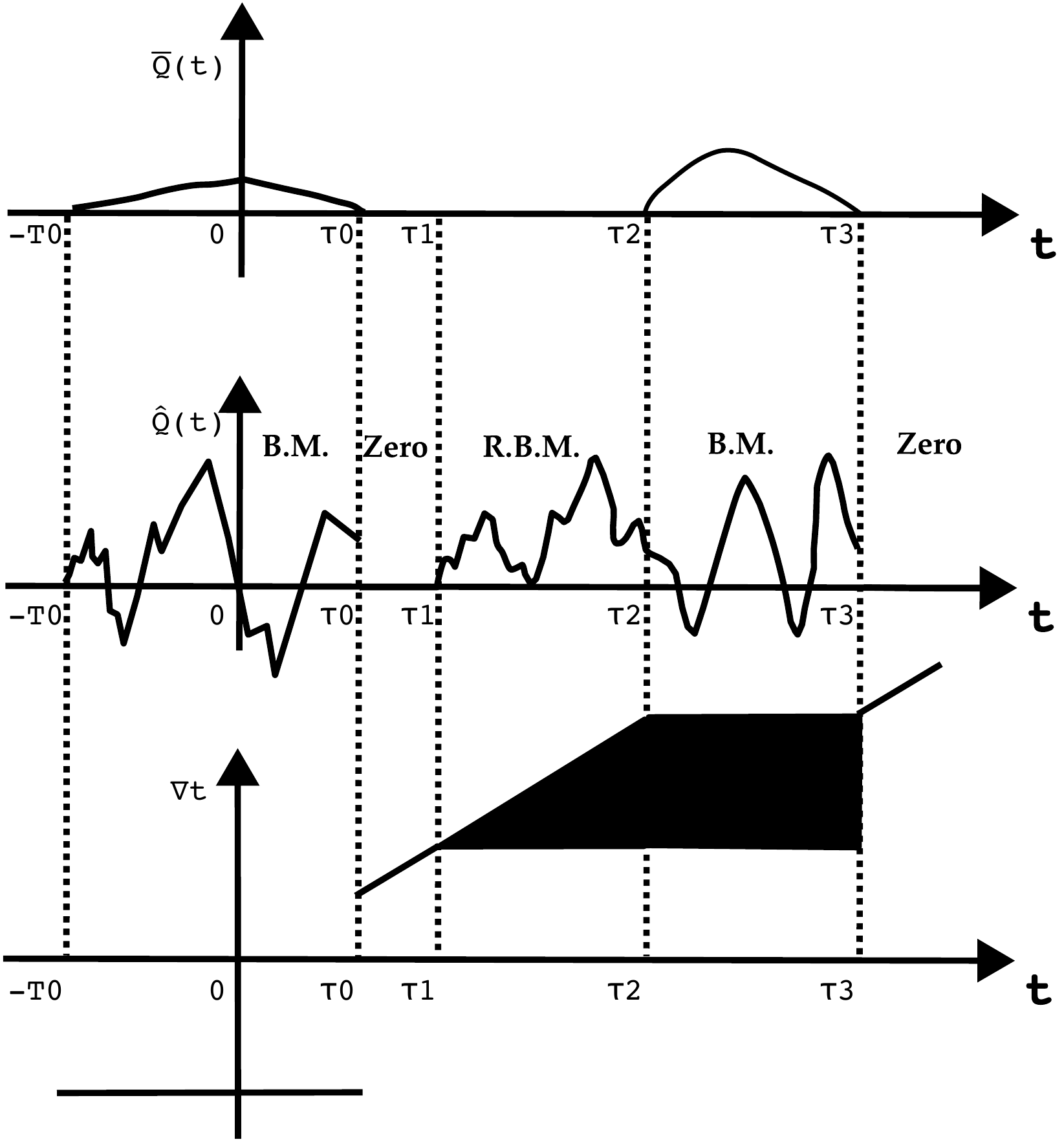}
\caption{An example of a $\D_{(i)}/GI/1$ queue that will undergo multiple
  ``regime changes''. The diffusion limit switches between a free Brownian motion (BM), a reflected Brownian motion (RBM), and the zero process.}
\label{fig:diffusion}
\end{figure}

\noindent 2. We noted in the remarks after Theorem
\ref{thm:queue-length-fluid} that the fluid limit can change between
being positive and zero in the arrival interval for a completely
general $F$. One can then expect the diffusion limit to change as
well, and switch between being a `free' diffusion, a reflected
diffusion and a zero process. This is indeed the case. Figure
\ref{fig:diffusion} illustrates this for the example in Figure
\ref{fig:fluid}. Note that $\forall t \in [-T_0,\t_1)$
$\Psi(\bar{X})(t) = -\bar{X}(-T_0)$, implying that the set
$\nabla_t^{\bar{X}}$ is a singleton. On the other hand, at $\t_1$
$\nabla_t^{\bar{X}} = \{-T_0,\t_1\}$. For $t \in (\t_1,\t_2]$,
$\Psi(\bar{X})(t) = 0 = \bar{X}(t)$, implying that $\nabla_t^{\bar{X}}
= (\t_1,t]$. On $(\t_2,\t_3)$, $\Psi(\bar{X})(t) = 0$, but $\bar{X}(t)
> 0$, so that $\nabla_t^{\bar{X}} = (\t_1,\t_2]$. Finally, the fluid
queue length becomes zero when the fluid service process exceeds the
fluid arrival process in $[\t_3,\infty)$, implying that
$\Psi(\bar{X})(t) = -(F(t) - \mu t) > 0$. It can be seen that
$\nabla_t^{\bar{X}} = \{t\}$ in this case. \vspace{3 pt}

\noindent 3. In \cite{HoJaWa12}, in fact, we prove joint
convergence of $(\hat Q^n,\tilde Y^n)$, in
the \emph{strong} $M_1$ topology, for the $\D_{(i)}/GI/1$ model.


\subsubsection{Why $M_1$, and not $J_1$?} \label{sec:M1J1}
We now discuss why we establish the diffusion limit in the space
$(\sD_{\lim}, M_1)$, and why it can't hold in the space $(\sD_{\lim},
J_1)$ in general. This section can be skipped on a first reading without any loss of continuity, though we encourage the reader to read it for a better understanding of Theorem \ref{thm:queue-length-diffusion}.

There are several equivalent definitions of \textit{convergence in the
  $M_1$ topology} (the interested reader is directed to \cite{Sk56,
  Wh01b, Wh01} for an in-depth study.) A simple characterization of convergence in $M_1$ for processes with range in $\mathbb{R}$ is the following involving the number of visits to a strip $[\alpha,\beta] \subset \mathbb{R}$ in an interval $[t_1,t_2] \subset [\eta,\infty)$. Let $y \in \sD$ (or $\sD_{\lim}$) and suppose there are $N+1$ points $t_1 \leq t_{(0)} < t_{(1)} < \ldots < t_{(N)} \leq t_2$ such that either $y(t_{(0)}) \leq \alpha, y(t_{(1)}) \geq \beta, y(t_{(2)}) \leq \alpha, \cdots$, or $y(t_{(0)}) \geq \beta, y(t_{(1)}) \leq \alpha, y(t_{(2)}) \geq \beta, \ldots$ Then, there are $N$ visits to the strip in $[t_1,t_2]$. Let $\nu_{[t_1,t_2]}^{[\alpha,\beta]}(y) \longmapsto \mathbb{N}$ be the number of visits to the strip $[\alpha,\beta]$ in $[t_1,t_2]$ by the function $y$. Definition \ref{lem:m1-conv} summarizes this characterization \cite{Wh01b}.

\begin{definition}[\emph{Convergence in $M_1$}] \label{lem:m1-conv}
~Let $y,y_n$ be elements of $\sD$ and $d_{M_1}(\cdot,\cdot)$  the $M_1$ metric. Then, $d_{M_1}(y_n,y) \longrightarrow 0$ as $n \rightarrow \infty$ if and only if
\[
\nu_{[t_1,t_2]}^{[\alpha,\beta]}(y_n) \longrightarrow
\nu_{[t_1,t_2]}^{[\alpha,\beta]}(y) \text{ as } n \rightarrow \infty.
\]
\end{definition}

Convergence in the $J_1$ topology, on the other hand, can be seen as a ``relaxation'' of the definition of convergence in the uniform metric topology. Specifically, let $z_n,z$ be elements of the space $\sD_{\lim}[\eta,\infty)$. Fix $T\in [\eta,\infty)$ that is a continuity point of $z$, and let $\|\cdot \|$ be the local uniform metric on the interval $[\eta,T]$. Define $\L$ to be the set of all non-decreasing continuous homeomorphisms from $[\eta,T]$ to itself. Then, convergence in $J_1$ can be defined as follows.
\begin{definition}[\emph{Convergence in $J_1$}] \label{def:j1-conv}
~There exists a sequence $\{\l_n\} \subseteq \L$ such that $\|\l_n - e \| \longrightarrow 0$ as $n \rightarrow \infty$, where $e$ is the identity map, $d_{J_1}(z_n,z) \longrightarrow 0$ as $n \rightarrow \infty$ if and only if $\| z_n \circ \l_n - z \circ e \| + \| \l_n - e \| \longrightarrow 0$ as $n \rightarrow \infty$.
\end{definition}
It is well known that the $M_1$ topology is weaker than the $U$ (uniform) or $J_1$ topologies, and processes converging in $M_1$ need not converge in $U$ or $J_1$.

As already stated, the diffusion limit for the queue length process is
obtained in the space $\sD_{\lim}$ when endowed with the $M_1$
topology because the directional derivative reflection mapping lemma
(Lemma \ref{lem:mama95}) that we use yields convergence in the $M_1$
topology alone. Intuitively, the reason the convergence result holds
only in $M_1$ is that asymptotically $y_n$ converges to a continuous
process, and it is well known that continuous processes can converge
to discontinuous limits only in the $M_1$ topology. To make this
intuition concrete, we give a counterexample that shows that convergence in $J_1$ is not possible in this case.

\begin{figure}[h]
  \centering
  \includegraphics[scale=1.2]{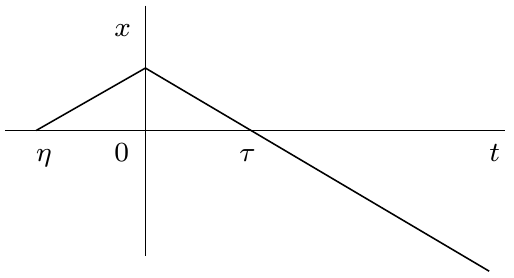}
  \caption{This $x \in C[\eta,\infty)$ corresponds to the fluid netput process, when $F$ is uniform.}
  \label{fig:x-ample}
\end{figure}

It will suffice to show that for some $\e > 0$ at least one of the terms in the expression $d_{J_1}(z_n,z) = \| z_n \circ \l_n - z \circ e \| + \| \l_n - e \|$ exceeds $\e$. Define the process $\ty_n$,
\[
\ty_n = \Psi(\sqrt{N_n} x + y) - \sqrt{N_n} \Psi(x),
\]
where $x$ is the function in Figure \ref{fig:x-ample}, and $y$ is a Brownian motion. We show that there is a non-empty set of points in the vicinity of $\t$ where the normed distance $d_{J_1}(\ty_n,\ty) > \e$, for any $\e >0$. Recall that $\ty = \sup_{s \in \nabla_{\cdot}^{x}} (-y(s))$.  The next proposition formalizes this argument.

\begin{proposition} [Non-convergence in $J_1$]
\label{thm:j1-convergence}
~Let $x$ be the function in Figure \ref{fig:x-ample}, $\{y_n\} \subset \sD_{\lim}[\eta,\infty)$ and $y \in \sC[\eta,\infty)$ is a Brownian motion, such that $y_n \stackrel{a.s.}{\to} y$ in $(\sD_{\lim}[\eta,\infty),U)$. Then, the process $\ty_n = \Psi(\sqrt{N_n} x + y_n) - \sqrt{N_n} \Psi(x)$ does not converge to $
\ty$ in the $J_1$ topology as $n \rightarrow \infty$.
\end{proposition}
The proof is available in the Appendix.

Thus, we see that the process $\ty_n$ does not converge to the
directional derivative of the reflection map in the $J_1$ topology
(and hence even the uniform topology), necessitating the use of the
$M_1$ topology. This result clearly implies that $\tY_n$ does not
necessarily converge to $\tY$ in the $J_1$ topology. Thus, we have a
situation where the limit process is discontinuous and the limit
result must be proved in the $M_1$ topology in full generality. 

\subsubsection{Busy Time Process} \label{sec:busy-time}
Recall from Section \ref{sec:fluid} that the fluid-scaled busy time
process $B^n(t)$ converges to a
continuous process $\bar{B}(t)$ as $n \to \infty$, in Corollary \ref{cor:busy-time-fluid}. Now, define the diffusion-scaled busy time process as
\begin{equation}
\label{busy-time-diffusion-scale}
\hat{B}^n \,:=\, \sqrt{N_n} (\bar{B} - B^n).
\end{equation}
Note that from the definitions of $B^n(t)$ and $\bar{B}(t)$ it follows that
$\hat{B}^n(t) = 0, ~ \forall t < 0$. As might be expected, the
diffusion refinement displays the same non-stationarity observed above.

\begin{corollary}
\label{cor:busy-time-diffusion}
~The diffusion scaled busy time process converges to a (directional
derivate) reflected Gaussian process as $n \to \infty$,
\[
\hat{B}^n \Rightarrow \hat{B} \,:=\, \frac{1}{\mu} \max_{s \in
  \nabla_{\cdot}^{\bar{X}}} (-\hat{X}(s)), ~\text{ in }~ (\sD_{\lim},M_1).
\]
\end{corollary}
\Proof
Recall that
\(
B^n(t) = s t \mathbf{1}_{\{t \geq 0\}} - I^n(t).
\)
~Substituting this and $\bar{B}$ from \eqref{busy-time-fluid} in the definition of $\hat{B}^n$, and rearranging the expression, we obtain
\(
\hat{B}^n = \frac{1}{\mu} \tilde{Y}.
\)
~A simple application of Theorem \ref{thm:queue-length-diffusion} then
provides the necessary conclusion.
\EndProof

Observe that $B^n(t)$ is approximated in distribution by $\hat{B}$ as
\(
B^n(t) \stackrel{d}{\approx} \bar{B}(t) - \frac{1}{\sqrt{N_n}} \hat{B}(t),
\)
~where $Z^n \stackrel{d}{\approx} Z$ is defined to be
\(
\bbP(Z^n \leq x) \approx \bbP(Z \leq x),
\)
and the approximation is rigorously supported by an appropriate weak
convergence result. Expressing the result in this manner exposes the
fact that the probability that the queue idled at any time up to $t >
0$ is
\[
\bbP(I^n(t) > 0) \approx \bbP(\hat B(t) > \sqrt N_n \bar I(t)).
\]
Note that the approximation is most accurate when the system is such
that it started with some workload at time $t = 0$. 

\section{Transitory Traffic Models} \label{sec:traffic-models}
In this section, we study three different traffic models for
transitory queueing systems, all of which satisfy Assumption \ref{assume:arrival-process}. We note that these models
are very different in nature and there could be many more that satisfy
the assumptions. Recall that the collection of arrival epochs $\mathbf
T_n := (T_{n,1}, T_{n,2}, \cdots, T_{n,N_n})$ is an instance of a
finite point process. We
assume, without loss of generality, that the support of the arrival
epochs is $[0,1]$ in this section. We give a brief description of the
models before continuing to a detailed description of each.

First, we model the arrival epochs as independent random
variables. Customers enter the queue in order of the sampled arrival
times, where each arrival time is sampled from a customer dependent
distribution. As noted before in Proposition
\ref{prop:empirical-process-limit} we studied a special case of this model in \cite{HoJaWa12},
where the arrival epochs were assumed to be i.i.d. We call this
the general $\D_{(i)}$ traffic model, for reasons elaborated on
below. See Section \ref{sec:delta-i-model}.

Classical queueing theory has focused extensively on modeling traffic
by renewal processes. In our second model, we assume that the joint distribution of the arrival
epochs is determined by conditioning the arrival epochs of a renewal
point process on an appropriate set of interest. We prove fluid and
diffusion limits for this model, and establish a close connection with
the asymptotics of the $\D_{(i)}$ model. See Section \ref{sec:conditioned-renewal}.

Finally, we consider a model of scheduled traffic with uncertainty,
where the realized arrival epoch is different
from the scheduled epoch, due to the fact that users sometimes arrive 
before or after the scheduled time. We model this variation between the realized and scheduled arrival times by a uniformly distributed random variable
with zero mean. We present fluid and diffusion limits to this model as
the population size tends to infinity. In particular, it is most interesting that this
model can be reduced to the general $\D_{(i)}$ model. See Section \ref{sec:scheduled}.

\subsection{The General $\D_{(i)}$ Model} \label{sec:delta-i-model}
The $\D_{(i)}$ traffic model assumes a product form for the joint
distribution of the arrival epochs. Independent arrival epochs are a
natural assumption to impose, especially while modeling a large number
of customers who independently take decisions on when to arrive at a
queue. In general, the marginals need not be identically distributed,
as individuals may have differing assessments on when to arrive at a
queue. We call this the \emph{general} $\D_{(i)}$ traffic
model. Customers enter the queue in order of the sampled arrival
times, so that the inter-arrival times are the difference of order
statistics, hence the term $\D_{(i)}$. In this section, we
comprehensively study the population acceleration (PA) fluid and diffusion limits for this model. 

Assume that $N_n$ is a deterministic, non-decreasing sequence of
natural numbers representing the population size in the $n$th system. Recall that $\mathbf T_n$ is an instance of the finite point
process, determining the arrival epochs. Without loss of generality of the domain, let $F_{n,i} : [0,1] \to [0,1]$,  for each $i \in \{1,
\ldots, N_n\}$, represent the arrival time distribution of 
customer $i$ in the $n$th system. Notice that $\{F_{n,i}, \, \, i  \in
\{1, \ldots, N_n\}\}$ $\forall n \in \bbN$ forms a triangular array of
distribution functions. As noted above the joint distribution of $\mathbf T_n$ is of product form, or formally
\(
\bbP(\mathbf T_n \in \Pi_{i=1}^{N_n} [0,t_i]) = \Pi_{i=1}^{N_n} F_{n,i}(t_i),
\)
for any Borel set $\Pi_{i=1}^{N_n} [0,t_i] \subset
[0,1]^{N_n}$. 

Intuitively, one might expect that any fluid limit to the arrival process \eqref{eq:arrival-process} would be an average over the individual sampling distributions. We first formalize this notion, by placing the following restriction on the $F_{n,i}$.  Let $\mathcal{K} :=
[0,1]$ represent the index set of customers, and $(\sK,\sB(\sK),m)$
represent the sample space of the indices, where $\sB(\sK)$ is the
Borel $\s$-algebra on $\sK$ and $m$ is the Lebesgue measure. Let
$\sL_{[0,1]}$ be the space of all distribution functions with support $[0,1]$. Following \cite{DuFr1967}, we define a \emph{random distribution function} as a mapping 
$\Upsilon: \sK \to \sL_{[0,1]}$. Thus, $(F_s (t) := \Upsilon(s)(t), ~ t \in [0,1])$, is the
distribution function of customer $s$. Clearly, $\Upsilon$ induces a
sample space $(\sL_{[0,1]}, \sB(\sL_{[0,1]}), \sP)$ where
$\sB(\sL_{[0,1]})$ is the Borel $\s$-algebra containing the weak-$^*$
topology on $\sL_{[0,1]}$ and $\sP = m \circ
\Upsilon^{-1}$ is the measure induced on the space $\sL_{[0,1]}$. 

The \textit{average distribution function} $\bar F$ is now well
defined in relation to $\sP$ as
\(
\bar F(t) := \int_{F \in \sL_{[0,1]}} F(t) d \sP (F) = \int_{[0,1]}
\Upsilon(s)(t) m(ds) = \int_0^1 F_s ds.
\)
  Notice that $\Upsilon$ is a measure-valued stochastic process with domain $[0,1]$
  and range $\sL_{[0,1]}$. It is useful to view $\Upsilon$ in the following sense: it represents
a summary of the beliefs of all the possible customers (in the
universe of these models) who may choose to arrive per the distribution
they choose from $\sL_{[0,1]}$. While the total order property of
$\sK$ plays no role in our description of the population of
customers, it is not unusual to expect that customers ``close'' to each
other, in the sense of the Euclidean norm on $\sK$, should have
similar beliefs. Thus, we impose the condition that $\Upsilon$
satisfies
\(
\| \Upsilon(\omega_1) - \Upsilon(\omega_2) \| \leq K |\omega_1 - \omega_2|,
\)
~for any $\o_1, \o_2 \in \sK$ and  $K < \infty$ is some given
constant. In particular, in the simplest case where $\Upsilon(\o) =
\d_{F} (\o)$ for all $\o \in \sK$ and some $F \in \sL_{[0,1]}$ (i.e.,
in an i.i.d. sampling model), this condition is satisfied
automatically. Note that this \emph{not} the usual sample path continuity of a stochastic process, but is instead a constraint on the variation of the sample path.

Recall that \eqref{eq:arrival-process} implies that the cumulative arrival process in the $n$th system is 
\(
A^n(t) := \sum_{i=1}^{N_n} \mathbf 1_{\{T_{n,i} \leq t\}}\,\, \forall
t \in [0,1].
\)
The fluid-scaled arrival process is simply 
\(
\bar A^n := \frac{A^n}{N_n},
\)
and the diffusion-scaled arrival process is 
\(
\hat A^n := \sqrt{N_n} \left( \bar A^n - \bar F_n \right).
\)
Customers enter the queue in the order
of sampled times. Thus, the inter-arrival times in the $\D_{(i)}$
arrival model are the differences of the ordered arrival times, so
that $\t_{(n,i)} = T_{(n,i)} - T_{(n,i-1)}$, where $T_{(n,0)} =
0$. 

Without loss of generality, assume customer $i$ corresponds to the
point $i/N_n \in (0,1]$. In order to establish large population fSLLN and fCLT's, we need some ``control'' on the
\emph{average distribution function} of row $n$ in the array, defined
as
\(
\bar F_n (t) := \frac{1}{N_n} \sum_{i=1}^{N_n} F_{n,i}(t), \,\, \forall t \in [0,1].
\) 
 We start by proving some useful properties of this average distribution function. The following lemma shows that $\bar F_n$ converges to $\bar F$ as $n
 \to \infty$. 

\begin{lemma} \label{lem:avg-dist-limit}
~There exists a distribution function $\bar F$ such that
\begin{equation} \label{eq:avg-dist}
\bar F_n (t) \to \bar F(t) := \int_{\sK} F_p(t) m(dp), 
\end{equation}
uniformly on $[0,1]$ as $n \to \infty$.
\end{lemma}

 A straightforward calculation shows that the covariance function of
 the arrival process in the $n$th row of the array is $K_n(s,t) := E[\hat
A^n(s) \hat A^n(t)] = \frac{1}{N_n} \sum_{i=1}^{N_n} F_{n,i}(s \wedge
t) - F_{n,i}(s) F_{n,i}(t)$. The following lemma shows that
$K_n$ has a well defined limit as $n \to \infty$.

\begin{lemma} \label{lem:avg-cov-limit}
~There exists a function $K(s,t)$ such that,
\begin{equation} \label{eq:cov-limit}
K_n(s,t) =  \frac{1}{N_n} \sum_{i=1}^{N_n} F_{n,i}(s \wedge
t) - F_{n,i}(s) F_{n,i}(t) \to K(s,t) := \int_{\sK} \left( F_p(s
  \wedge t) - F_p(s) F_p(t) \right) m(dp),
\end{equation}
as $n \to \infty$, uniformly for all $s,t \in [0,1]$.
\end{lemma} 
The proofs of Lemma \ref{lem:avg-dist-limit} and \ref{lem:avg-cov-limit}
are in the appendix.

The fSLLN theorem is a generalization of the Glivenko-Cantelli Theorem \cite{Du10}
to triangular arrays of non-identically distributed random variables. In Theorem \ref{thm:arrival-slln} we show that
the normalized arrival process $\bar A^n$ converges to $\bar F$ in $\sD_{\lim}$, uniformly on
compact sets of $[0,\infty)$. We prove this result by demonstrating the uniform convergence of the
sample paths of the empirical distribution. We make the reasonable
assumption that none of the distribution functions $\bar F_n$ share
discontinuity points in the support. This implies that the limit $\bar
F$ is (almost surely) continuous, allowing us to prove convergence in
the uniform metric. 

\begin{theorem}[Glivenko-Cantelli Theorem for Triangular Arrays] \label{thm:arrival-slln}
~The fluid scaled arrival process $\bar A^n = \frac{A^n}{N_n}$
satisfies a functional strong law of large numbers,
\[
\bar{A^n} \stackrel{a.s.}{\rightarrow} \bar F \text{ in } (\sD_{\lim}, U),
\]
as $n \to \infty$. 
\end{theorem}
The proof is presented in the appendix.

\noindent \textbf{Remarks.} Versions of this theorem have been proved in the literature and we draw attention, in particular, to
Theorem 1 of \cite{We1981}. There it was shown that $\bar A^n$ and
$\bar F_n$ converge to the same limit point in the Prokhorov metric (\cite{Bi68,Wh01}) on
$\sL_{[0,1]}$. However, this result does not explicitly identify the
fluid limit process. Our construction of the empirical distribution via random distribution
functions allows us to do this. 

A requirement in the proof of Theorem \ref{thm:arrival-slln} is that the sequence
$\{N_n, n\geq 1\}$ must satisfy $\sum_{n=1}^{\infty} \frac{1}{N_n^2} <
\infty$, implying that $N_n = O(n^{1+\d})$ for $\d \geq 0$, and this will
play a role in the proof of the fCLT. 

As a consequence of our construction of the empirical distribution function 
space and Lemma \ref{lem:avg-dist-limit}, we can explicitly identify the
Gaussian limit process in our setting. We first present the fCLT, and
then elaborate on the diffusion limit process.

\begin{theorem} [Empirical Process Limit for Triangular Arrays]\label{thm:arrival-fclt}
~The centered arrival process $\hat A^n$ satisfies a functional central
limit theorem,
\[
\hat A^n \Rightarrow \tilde W \text{ in } (\sD_{\lim},U),
\]
as $n \to \infty$, where $\tilde W$ is a mean zero Gaussian process with covariance
function $K(s,t)$ defined in \eqref{eq:cov-limit} and continuous
sample paths.
\end{theorem}
The proof can be found in the appendix.

\noindent \textbf{Remarks.} Our result is a generalization of Hahn's Central Limit Theorem
(Theorem 2 in \cite{Ha1978}) to nonidentically distributed random
elements of $\sD_{\lim}$. We also draw attention to Theorem
1.1 of \cite{Sh1979} that proves the existence of an empirical process
limit for triangular arrays, under the sufficient condition that an
appropriate covariance function exists. However, it does not specifically
identify the Gaussian process limit, something that is crucial for
this paper. Once again, our identification of the empirical distribution
function space with random distribution functions enables this identification.

The covariance structure of the process $\tilde W$ is
interesting in itself, and we make the following
observations. First, notice that the covariance function is an
\emph{average} of the covariance functions of the Brownian Bridge
processes $W^0 \circ F_p$ (with $p \in \sK$) where $W^0$ is a
standard Brownian Bridge. In a sense, these are Brownian Bridge processes associated with empirical processes of random samples from the function $F_p$. Second, differentiating the expression for
$K(t,t)$ with respect to $t$ we have
\(
\frac{d K(t,t)}{dt} = \int_{\sK} \left( f_p(t) - 2 f_p(t) F_p(t)
\right) m(dp),
\)
 where $f_p$ is the density (or at least the right-derivative) of the
 distribution function $F_p$. This is the average of the infinitesimal variance
  of the Brownian Bridges $W^0 \circ F_p$. Recall that the
  infinitesimal mean and variance of a diffusion process are defined
  as $\frac{E[X(t+h) - X(t) \vert X(t) = x]}{h} \to \mu(t,x)$ and
  $\frac{E[|X(t+h) - X(t)|^2 \vert X(t) = x]}{h} \to \sigma^2(t,x)$ as
  $h \to 0$ (resp.). For the Brownian Bridge process $W^0 \circ F_p$, it is well known that
  the infinitesimal mean and variance are (for a fixed $p \in \sK$) 
\begin{eqnarray*}
\m_p(t,y) &=& \frac{-y f_p(t)}{1 - F_p(t)} \\
\s^2_p(t,y) &=& f_p(t).
\end{eqnarray*}

Further, it can be shown that the mean and variance of the Brownian
Bridge satisfies the following o.d.e's:
\begin{eqnarray*}
\frac{d}{dt} \bbE \left[ W^0\circ F_p(t) \right] &=& \bbE\left[ \m_p(t,W^0\circ F_p
  (t))\right ] = \frac{-f_p(t)}{1 - F_p(t)} \bbE[W^0 \circ F_p(t)] = 0\\
\frac{d}{dt} \text{Var} \left(W^0 \circ F_p(t) \right) &=& \bbE \left[ \sigma^2(t,W^o
  \circ F_p(t)) \right ] + 2 \bbE\left[ W^0 \circ F_p(t) \times \m_p(t,W^0\circ
F_p(t)) \right]\\ &=& f_p(t) - 2 f_p(t) F_p(t).
\end{eqnarray*}

Comparing the variance derivative above with $\frac{d K(t,t)}{d t}$,
we conjecture that the $\tilde W$ is a Gaussian diffusion process
with infinitesimal generator equal to the average of the infinitesimal
generators of the Brownian Bridges $W^0 \circ F_p$. However, we have not been able to verify that the process is Markov with respect to its natural filtration to make a definitive conclusion.

A particular case of interest is when the $\{T_{n,i}\}$ are
i.i.d. drawn from a common continuous distribution $F \in \sL_{[0,1]}$. This
result, of course, is the standard fCLT for the empirical process (see
\cite{Bi68,WeSh2009,Wh01} for a deeper exposition).

\begin{corollary}\label{cor:iid-sample}
~For each $n \geq 1$, let $\{T_{n,i},\, i = 1,\ldots,N_n\}$ be a
triangular array of i.i.d. random
samples drawn from a distribution $F$. Then, as
$n \to \infty$
\[
\hat A^n \Rightarrow W^0\circ F ~~\text{in} (\sD,U).
\]
Here $W^0$ is the standard Brownian Bridge process defined on the common sample space.
\end{corollary}

A formal proof of this result is standard and omitted (see Chapter 13
of \cite{Bi68}). However, it is also
straightforward to see this from Theorem \ref{thm:arrival-fclt} by setting
$F_p = F$ for all $p \in [0,1]$. It can be readily verified that the Gaussian
process $\tilde W$ is equal in distribution to a Brownian Bridge
process.

\subsection{Conditioned Renewal Model} \label{sec:conditioned-renewal}
The most common traffic model assumed in the queueing theory
literature is the renewal model. In this section, we
consider a model of traffic for transitory queueing systems that is
related to renewal traffic models. Specifically, we allow the arrival
process to be a renewal process conditioned on the even that $N_n$
arrivals occur in some finite time horizon. First, we recall that there is a strong
connection between the i.i.d. $\D_{(i)}$ arrival process and the conditioned Poisson
renewal process. Next, we show that even though this property is not
satisfied by renewal processes in general, asymptotically we obtain
the same fSLLN and fCLT limit processes as the population size scales
to infinity. This appears to be a new result and should be of wider
interest. Without loss of generality we fix $N_n = n$ in this section.
\subsubsection{Relation Between Conditioned Poisson and $\D_{(i)}$
  Models}
Consider a renewal point process $(M(t), t\geq 0)$ defined with respect to $(\Omega, \sF, \bbP)$. Let
$(\lambda(t), t \geq 0)$ be an integrable, non-negative function
defined to be the arrival rate of $M$. Therefore, $\G(t) := \int_0^t \l(s) ds$ is
the mean cumulative arrival process. We first note the following.

\begin{lemma} \label{lem:dist-func-from-poisson}
~Let $\G : [0,\infty) \to [0,\infty)$ be the mean cumulative arrival
process of a Poisson process. Then, for a fixed $T > 0$,
\begin{equation} \label{eq:dist-func-from-poisson}
F(t) := \frac{\G(t)}{\G(T)} ~~\forall t \in [0,T],~~
\end{equation}
is a continuous probability distribution function.
\end{lemma}
It is straightforward to verify this result and we omit a proof. The ordered statistics
(\textsc{OS}) property of point processes provides the connection
between the i.i.d. $\D_{(i)}$ and Poisson processes.

\begin{definition}[Property \textsc{OS}] \label{def:OS}
~Conditioned on
$\{M(T) = n\}$, the event epochs $(T_1,
\ldots, T_n)$, are distributed as the ordered statistics of $n$
independent and identically distributed random variables with
distribution $F(t)$, for $t \in [0,T]$.
\end{definition}

By Theorem 1 of \cite{Li1985}, $M$ possesses the
\textsc{OS} property if and only if it is a Poisson process (see
\cite{Fe79} as well). Notice that this distributional relationship is
true for every $n \geq 1$. By Kolmogorov's Extension Theorem, there
exists a stochastic process $\hat M^n(t)$ such that for any
partition $0 < t_1 < \cdots < t_d < T$ and $(x_1,\cdots, x_d) \in
\bbR^d$, $\bbP(\hat M^n(t_1) \leq x_1, \ldots, \hat M^n(t_d) \leq x_d)$
\begin{equation*}
=\bbP \left(\frac{1}{\sn}(M(t_1) - n F(t_1)) \leq x_1,
  \ldots,\frac{1}{\sn}(M(t_d) - n F(t_d) ) \leq x_d | M(T) = n \right).
\end{equation*}
The \textsc{OS} property implies that we
can easily obtain an fCLT for the conditioned Poisson process.

\begin{theorem} \label{thm:conditioned-poisson-fclt}
~The sequence of processes $\{\hat M^n \}, ~ n \geq 1,$ satisfies a functional central limit theorem,
\[
\hat{M}_n \Rightarrow W^0 \circ F \,\, \text{ in } (\sD, U),
\]
as $n \to \infty$, where $W^0$ is a standard Brownian Bridge process defined on the same
sample space as $M$
\end{theorem}
The proof is simple and in the appendix. The implication of this
theorem is that the conditioned Poisson and $\D_{(i)}$ traffic models
are equivalent in distribution. In fact, verifying that an observed
traffic sequence satisfies the \textsc{OS} property is
sufficient to conclude that the arrival process is Poisson. A thorough
study of statistical tests for this purpose is presented in \cite{KiWh2013a}.\\

\noindent \textbf{Remarks} 1. It is important to note that this limit result fundamentally differs from the
standard diffusion limit for non-homogeneous Poisson processes, which we review here. Let $N(\cdot)$ be a
unit rate Poisson process. Then, $M(t) \stackrel{d}{=} N(\int_0^t \l(s) ds)$ is a non-homogeneous Poisson process. The diffusion approximation to
this process is developed by scaling the compensated
(Martingale) process $\hat N(t) := M(t) - \int_0^t \l(s) ds$ in an
appropriate manner. The commonly accepted approach is the
\emph{uniform acceleration} method developed in \cite{MaMa95}, where
the rate function $\l(s)$ is scaled by a constant $\e > 0$, so that
we obtain the scaled process
\(
\hat N^{\e} (t) := N\left( \frac{1}{\e} \int_0^t \l(s) ds\right) - \frac{1}{\e} \int_0^t \l(s) ds.
\)
~In \cite{MaMa95}, the Strong Approximation Theorem (see \cite{KoMaTu76})  is used to
prove that the sample paths of $\hat N^{\e}$ converges to those of a
standard Brownian motion as $\e \to 0$. Recall that the strong
approximation theorem implies that 
\(
\hat N^{\e}(t) = \frac{1}{\sqrt{\e}}W(t) + o\left( \frac{1}{\sqrt{\e}}
\right) ~~ a.s.
\)
~as $\e \to 0$. When $\l(s) = \l$ for all $s
\geq 0$, and $\e = 1/n$ for $n \in \bbN$, it is straightforward to see
that the standard Poisson process diffusion approximation is a special
case of this approach (see also \cite{Gl1998} for an overview of
strong approximation methods applied to queueing theory). 

Now, contrast this limit with Theorem
\ref{thm:conditioned-poisson-fclt}. Note that, by definition, $\hat
M^n$ is \emph{not} equivalent to the compensated Martingale process
$\hat N^\e$ (when $\e = 1/n$), as it is defined with respect to the
conditioned measure on the set $\{M(T) = n\}$, and not the full
measure $\bbP$. 
Furthermore, the limit can only hold in the weak sense, as the strong approximation only applies
to processes with independent increments. This is not a condition satisfied by $\hat M^n$. It
appears that obtaining a strong approximation (or rate of convergence)
result for the conditioned process is an open problem, and of independent interest.

\subsubsection{Functional Limit Theorems of Conditioned Renewal
  Processes}
Theorem 1 of \cite{Li1985} clearly shows that a non-Poisson renewal
process does not satisfy the \textsc{OS} property. However, in this
subsection we prove that the conditioned renewal process in fact converges to
a Brownian Bridge process when scaled appropriately. For simplicity,
we assume that the renewal process is time-homogeneous, but the
results extend easily to the general case. Without loss of generality
we also assume that $\l(t) = 1$ for all $t \geq 0$. First, we recall the
definition of a \textit{finitely exchangeable} sequence.

\begin{definition}
~Let $\{X_1, \ldots, X_n\}$ be a collection of random variables defined
with respect to the sample space $(\Omega, \sF, \bbP)$. Then, this
collection is said to be finitely exchangeable if
\(
\{X_1, \cdots, X_n\} \stackrel{D}{=} \{X_{\pi(1)}, \cdots, X_{\pi(n)}\},
\)
where $\pi : \{1, \ldots, n\} \to \{1, \ldots, n\}$ is a permutation
function on the index of the collection.
\end{definition}

Renewal processes satisfy the \emph{exchangable} (or \textsc{E})
property, as summarized in the following proposition.

\begin{proposition}[Property \textsc{E}] \label{prop:exchangeinter}
~Let $\xi_{i} : \Omega \to \bbR_+$ $i \in \bbN$ be a sequence of
i.i.d. positive random variables defined with respect to $(\Omega,
\sF, \bbP)$, such that $M(t) := \sup \{k > 0 |
\sum_{l=1}^k \xi_{l} \leq t\}, $ for all $t > 0$ is the associated
renewal counting process. Then, the finite collection $\Xi_{n} :=
(\xi_{1}, \ldots, \xi_{n})$ is \textit{finitely exchangeable} under the measure
conditioned on the event $\{M(T) = n\}$, for $T < \infty$ fixed.
\end{proposition}

A proof of this fact is in the appendix. Notice that finitely
exchangeable random variables are not \emph{infinitely} exchangeable and 
important results such as de Finetti's Theorem are unavailable. To
prove the functional limit theorems for the counting processes, we will first prove that certain
scaled partial sums of the exchangeable inter-arrival times
(conditioned) have Gaussian process weak limits. However, in contrast to classical functional
central limit theorem results, the conditioning increases the
complexity of the problem significantly, since for each $n$ the random
variables exist on different (but related) probability sample
spaces. The result for the counting process will follow by a
random time-change argument. \\

\noindent \textbf{Functional Strong Law of Large Numbers}
Consider a triangular array of random variables
$\Xi_n := (\xi_{n,i},~~i=1,\ldots,n)$ and $n \geq 1$, defined as the inter-arrival times of the
renewal events associated with a sequence of independent and
indistinguishable renewal processes, $\{M_n,~n \geq 1\}$. By Proposition
\ref{prop:exchangeinter}, we know that $\Xi_n$ is an exchangeable
array of random variables conditioned on the event $\{M_n(T) =
n\}$. The limit results are proved with respect to a conditional measure $\bar
\bbP$ that we construct in Section \ref{sec:lemmata} in the Appendix. This can be
skipped on a first reading by accepting the premise that such a
measure exists. In the ensuing, any reference to $\Xi_n$ is to be interpreted
with respect to the conditional measure $\bar \bbP$.

Let $\m_n := E_{\bar \bbP}[\xi_{n,i}] = E[\xi_{n,i} | M_n(T) = n]$ be the conditioned mean
of the inter-arrival periods; the exchangeable property implies that
these random variables are identically distributed. Our first result is a functional strong law
for partial sums of these random variables.

\begin{theorem} \label{thm:conditioned-fslln}
~Let $\bar S_n(t) := \sum_{l=1}^{\lfloor nt \rfloor} \xi_{n,l} \in
\sD_{\lim}$, $\forall t \in [0,1]$. Then, 
\[
~\bar S_n  \stackrel{\bar \bbP - a.s.}{\to} e ~~\text{in}~~ (\sD_{\lim},U)
\]
as $n \to \infty$, where $e:[0,1] \to [0,1]$ is the identity function.
\end{theorem}

This is an intuitively satisfying result, that provides strong evidence
that an fCLT along the lines of Theorem \ref{thm:conditioned-poisson-fclt}
is satisfied by a conditioned renewal process.\\

\noindent \textbf{Functional Central Limit Theorem}
Consider the standardized random
variables, $\{\phi_{n,l}, ~~ l = 1,\ldots,n\}$ defined with respect to $\Xi_n$:
\[
\phi_{n,l} := \frac{\xi_{n,l} - \mu_n}{\sqrt{n}}.
\]
The next theorem characterizes the sequence $\phi_{n,l}$ and shows
that the partial sums of these random variables satisfy a functional
central limit theorem. 


\begin{theorem} \label{thm:conditioned-fclt} 
~Let $\{\phi_{n,l}, ~~ l=1,\ldots,n\},~n \geq 1,$  be the triangular
  array of random variables defined above and $\hat M_n(t) :=
\sum_{i=1}^{\lceil n t \rceil} \phi_{n,i} ~~ \in
\sD_{\lim}~~\text{and}~~ \forall t \in [0,1]$. Then, the random
variables $(\phi_{n,1}, \ldots, \phi_{n,n})$ are exchangeable and satisfy:\\
\noindent (i) $\sum_{l=1}^n \phi_{n,l} \stackrel{\bar \bbP}{\to} 0$, \\
\noindent (ii) $\max_{1 \leq l \leq n} |\phi_{n,l}| \stackrel{\bar \bbP}{\to} 0$, \\
\noindent (iii) $\sum_{l=1}^n \phi_{n,l}^2 \stackrel{\bar \bbP}{\to} 1$, and\\
\noindent (iv) 
\(
\hat M_n \Rightarrow W^0 ~~\text{in}~~ (\sD_{\lim},U),
 \)
 as $n \to \infty$, where $W^0$ is a standard Brownian Bridge process.
\end{theorem}
\Proof
The exchangeability of $\phi_{n,i}$ follows directly by the fact that
$\xi_{n,i}$ is exchangeable. (i), (ii) and (iii) are proved  in Proposition \ref{prop:verification}
in the appendix. Then, by Theorem 24.2 of \cite{Bi68} $\hat M_n
\Rightarrow W^0$ in $(\sD,U)$. The extension to $(\sD_{\lim},U)$
follows from Lemma \ref{lem:d-lim}.
\EndProof

The conditions in Theorem \ref{thm:conditioned-fclt} are natural in the context of the
conditioned limit result we seek. Note that the conditioned limit result is akin
to proving a diffusion limit for a tied-down random walk (see
\cite{Li1970,Me1989}). The first condition here enforces a type of
``asymptotic tied down'' property. The second condition is a necessary and
sufficient condition for the limit process to be infinitely divisible
(see \cite{ChTe1958} for more on this). The third condition is
necessary to ensure that the Gaussian limit, when $t = 1$, has
variance $1$. Similar conditions have been observed to be
sufficient to prove central limit theorems for dependent random
variables (see, in particular, \cite{Mc1974,We1980}). 

The final step in describing the limit behavior of the conditioned
renewal process is to obtain a result that parallels Theorem
\ref{thm:conditioned-poisson-fclt}, and show that the ``counting''
counterpart of the partial sum process also converges to a Brownian Bridge.
To that end, we start with a definition. 

\begin{definition}[Counting Process]
~$L_n(t) := \sup \{ 0 \leq m \leq n | \bar S_n(\frac{m}{n}) := \sum_{l=1}^{m} \xi_{n,l} \leq t\}$ is the
standard counting counterpart to the partial sum process. 
\end{definition}

Now, the main theorem of this section proves that the counting
process satisfies Assumption \ref{assume:arrival-process}.

 \begin{theorem} \label{thm:counting}
~(i) $\bar L_n := \frac{L_n}{n} \stackrel{a.s.}{\to} e ~~ \text{in} ~~
(\sD_{\lim},U)$ as $n \to \infty$, where $e : [0,\infty) \to [0,\infty)$ is the identity map.\\
~ (ii) $\sn \left( \bar L_n - e \right) \Rightarrow -W^0 ~~ \text{in}
~~ (\sD_{\lim},U)$ as $n \to \infty$, where $W^0$ is the Brownian
Bridge limit process observed in Theorem \ref{thm:conditioned-fclt}.
 \end{theorem}

Part (i) shows that the conditioned renewal traffic model converges to the uniform distribution function
on $[0,1]$, in a large population limit. Part (ii), in turn, proves
that the diffusion scaled conditioned renewal traffic model satisfies
Assumption (b) in Assumption \ref{assume:arrival-process}. The proof is a
consequence of the random time change theorem (see chapter 17 of \cite{Bi68}), and relegated to the appendix. This is an intuitively satisfying result as we know that the Poisson
renewal model certainly satisfies the same limit. Further,
as a result of Theorem \ref{thm:counting}, it is obvious that the
large population approximations to the performance metrics of a ``conditioned renewal'' transitory queueing model is asymptotically
equal in distribution to an equivalent $\D_{(i)}$ transitory queueing
model. 

\subsection{Scheduled Arrivals with Epoch Uncertainty} \label{sec:scheduled}
Traffic scheduled to arrive at regular intervals is a common occurrence
in many service systems. Often, schedules are made for a finite period
of time, and the traffic pattern is transitory in nature. For example, 
hospital outpatient units schedule patients at particular times during the day (typically 8AM to 8PM). Another classic example of scheduled
traffic is air traffic arrivals. However, while the arrivals may be scheduled, it is often
the case that there is some randomness in the realized arrival time: users can arrive a little before or after the scheduled
arrival time. The earliest description of a model of such arrival behavior appears in
\cite{CoSm1961} where it was introduced as ``a regular arrival process
with unpunctuality''. Recent work in \cite{ArGl2012} studied this
model with heavy tailed uncertainty and demonstrated convergence to a
fractional Brownian motion with Hurst index $< 1/2$. In this section,
we present a novel and intuitive model of scheduled arrivals with
uncertainty on a finite interval, and demonstrate its connection with
the $\D_{(i)}$ traffic model.

For simplicity, let the population size be $N_n = n \in \bbN$, and without loss of
generality we assume that arrivals take place over the interval $[0,1]$ at equal intervals. For
simplicity we assume the first arrival is scheduled at time 0, and the
last one at time 1. The $j$th user arrives at time $\t_{n,j} :=
j/n$; for simplicity, assume $0=\t_{n,1} \leq \t_{n,2} \cdots \leq
\t_{n,n} = 1 $. Let $\xi_{n,i}$ be a random variable uniformly distributed on
the interval $[-T,T]$, where $T$ is a constant to be defined. Then, the realized arrival time of user $j$ is
modeled as $T_{n,j} := \t_{n,j} + \xi_{n,j}$. Users can potentially enter the service
system in the interval $\sT := [-T,1+T]$, and the cumulative
number of arrivals by time $t \in \sT$ is $A^n(t) := \sum_{i=1}^n
\mathbf 1_{\{T_{n,i} \leq t\}}$. 

We now argue that this arrival process satisfies Assumption
\ref{assume:arrival-process}. Analogous to Section
\ref{sec:delta-i-model}, consider the average distribution function
\[
\bar F_n (t) := \frac{1}{n} \bbE[A^n(t)] = \frac{1}{n} \sum_{i=1}^n F(t - \t_{n,i}).
\]
Note that the summands are not the same distribution function, as the mean of
$T_{n,i}$ is $\t_{n,i}$. This is analogous to the definition of the
average distribution function in Lemma
\ref{lem:avg-dist-limit}. Our first result argues that there exists a functional limit to $\bar F_n$ as $n
\to \infty$.

\begin{proposition} \label{prop:scheduled-avg-dist-limit}
~Let $\bar F_n$ be the average distribution function, for a given $n
\geq 1$. Then, for a fixed  $T \in [0,0.5]$ and $t \in \sT$ 
\[
\bar F_n(t) \stackrel{a.s.}{\to} \bar F(t) := \begin{cases}
\frac{(t+T)^2}{4 T}, & ~~ -T \leq t \leq T,\\
	t, & ~~ T < t \leq 1 -T, \\
	\begin{split} &\frac{t+T}{2 T} - \frac{t^2 - T^2}{2 T}\\ &+
          \frac{(t - T)^2}{4 T} - \frac{1}{4 T} + (t-T), \end{split} & ~~ 1-T < t \leq 1+T,
\end{cases}
\]
uniformly on $[0,1]$ as $n \to \infty$.
\end{proposition}

\begin{figure}
\centering
\includegraphics[scale=0.3]{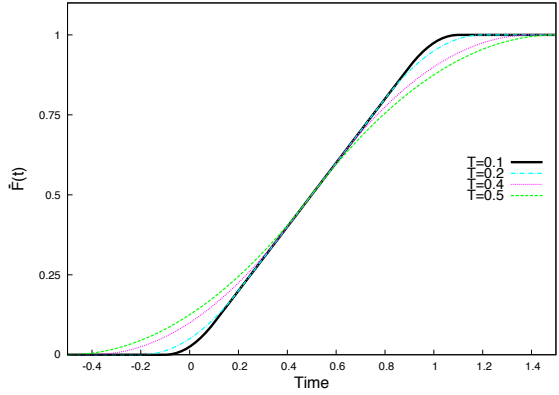}
\caption{The ``population'' average arrival distribution function for
  different values of $T$.}
\label{fig:scheduled}
\end{figure}

The proof is available in the appendix. Figure \ref{fig:scheduled}
depicts $\bar F$ for different values of $T$. Notice that the support of
the population mean distribution function depends on the value of $T$, and the larger the
value of $T$, the earlier and later arrivals can occur to the system (obviously). Interestingly enough, one
can also show that the limit population (or \emph{mean field})
distribution is an average over the individual distribution functions
for each user. Following Section \ref{sec:delta-i-model}, let $\sK :=
[0,1]$ represent the universe of all possible users to the queueing
system. Then, for $p \in \sK$, we know that
\begin{equation} \label{eq:uncert-dist}
F_p(t) := \frac{t-p+T}{2T} ~~\text{for}~~ p-T \leq t \leq p+T,
\end{equation}
is the arrival distribution function of customer $p$. That is, customer $p$ arrives at \emph{time}
$p$, with a uniform uncertainty distribution centered at $p$. The
following corollary shows that this population average coincides with
the distribution $\bar F$ in Proposition
\ref{prop:scheduled-avg-dist-limit}, when $T \in [0,0.5]$. The proof is a simple integration
argument, and is relegated to the appendix.

\begin{corollary} \label{cor:scheduled-avg-dist-limit}
~Let $F_p$ defined in \eqref{eq:uncert-dist} be the arrival distribution
associated with user $p \in \sK$. Then,
\[
\int_0^1 F_p(t) m(dp) = \begin{cases} 
\frac{(t+T)^2}{4 T}, & ~~ -T \leq t \leq T,\\
	t, & ~~ T < t \leq 1 -T, \\
	\begin{split} &\frac{t+T}{2 T} - \frac{t^2 - T^2}{2 T}\\ &+
          \frac{(t - T)^2}{4 T} - \frac{1}{4 T} + (t-T), \end{split} & ~~ 1-T < t \leq 1+T,
\end{cases}
\]
where $m(\cdot)$ is the Lebesgue measure on the set $\sK$.
\end{corollary}

This is precisely the condition that needs to be satisfied for the
generalized Glivenko-Cantelli result to be true in Theorem
\ref{thm:arrival-slln}. Let $F_{n,i}(t) := F(t - \t_{n,i})$ and extend
the support of $F_{n,i}(t)$ to the interval $[-T,1+T]$ such that
$F_{n,i}$ is zero outside the interval
$[-T+\frac{i}{N_n},T+\frac{i}{N_n}]$. The
scheduled arrival model is, therefore, a special case of the general $\D_{(i)}$
traffic model. We claim the following theorem as an immediate
consequence of Theorem \ref{thm:arrival-slln}.

\begin{theorem} \label{thm:scheduled-slln}
The fluid-scaled arrival process satisfies a functional strong law of
large numbers:
\[
\bar{A^n} \stackrel{a.s.}{\rightarrow} \bar F \text{ in } (\sD_{\lim}, U),
\]
as $n \to \infty$.
\end{theorem}

Similar arguments as above we can also prove that the sample
covariance function also converges to the average covariance, as
obtained in Lemma \ref{lem:avg-cov-limit}. Since the arguments are
similar, we skip the proof and present the final result. The
diffusion-scaled arrival process is denoted
\(
\hat A^n(t) := \sn \left( \frac{A^n(t)}{n} - \frac{1}{n} \sum_{i=1}^n
  \mathbf F( t - \t_{n,i}) \right).
\)

\begin{theorem} \label{thm:scheduled-fclt}
~The centered arrival process $\hat A^n$ satisfies a functional central
limit theorem,
\[
\hat A^n \Rightarrow \tilde W \text{ in } (\sD_{lim},U),
\]
as $n \to \infty$, where $\tilde W \in \sC$ is a zero mean Gaussian process with covariance
function $K(s,t)$, as obtained in Lemma \ref{lem:avg-cov-limit}.
\end{theorem}

As a final note, observe that there is an important distinction between the scheduled arrival 
and the $\D_{(i)}$ pre-limit traffic models. In the latter, the realized arrival
times are the ordered statistics of the sampled arrival times, while in the former this is not the
case. However, the natural (partial) ordering of the real numbers is
all that is required to establish the functional limits, and in the
limit as $n \to \infty$ any difference between these models is
``washed out''. The results in this and the previous section strongly
indicate that, in some sense, the $\D_{(i)}$ traffic model is
canonical to the study of transitory queueing systems. In the next
section we focus on a deeper study of sample path approximations
suggested for such models.

\section{Conclusions and Future Work}\label{sec:conclusion}

In this paper, we introduced a framework for \textit{transitory queueing systems}. We define transitory queueing models as those whose arrival process is a finite point process that satisfies (i) the empirical arrival process satisifes a fSLLN and the limit is a well defined probability distribution function, and (ii) the empirical arrival process satisfies a fCLT, converging to a tied down Gaussian process when appropriately centered and scaled.  Our attempt is to capture queueing scenarios where the queues are `transitory' either because the queue operates only for a finite time, or because there is only a finite population of customers that will arrive. Such scenarios are very difficult to study via classical queueing analysis, and in some cases even non-stationary Markovian models may not be appropriate. 

We first develop fluid and diffusion approximations to the system performance metrics as the population size increases, in a very general setting. More precisely, we only assume that the traffic model satisfies a fluid and diffusion limit, and that the service process is renewal. Under these conditions, we derive a fluid limit to the queue length process and show that it can switch between various regimes. We then derive the diffusion limit of the queue length process in terms of a directional derivative of the Skorokhod reflected fluid ‘netput’ process. This diffusion limit process switches between a free diffusion, a reflected diffusion and the zero process. The weak convergence proof is shown to hold in the $M_1$ topology on the space $\sD_{\lim}$. A novel feature of the diffusion limit is that it is a function of a Gaussian bridge process and a Brownian motion process. This appears to be a unique observation as such diffusion limits do not appear in the ``conventional heavy-traffic'' approximations to multi-server queues. 

Then, we introduce three natural `transitory' traffic models that satisfy the assumptions. In the $\D_{(i)}$ model, a finite population of customers choose (sample) a time of arrival at the queue in an independent manner from distribution functions that are potentially unique to each of them. The customers then enter in order of the sampled arrival times - thus the arrival epochs are ordered statistics. We developed generalizations of the Glivenko-Cantelli and Donsker's Theorems for triangular arrays, and show that the limit process has an intuitive and interesting interpretation. 

Next, we investigated the connection between classical queueing models and the $\D_{(i)}/GI/s$ queue. We rigorously show that a conditioned renewal process (appropriately scaled) converges to a Brownian Bridge process \emph{a la} the i.i.d. sampling $\D_{(i)}$ process. This implies that the performance metrics of a (conditioned) $GI/G/s$ queueing model are asymptotically equivalent in distribution to those of the $\D_{(i)}/GI/s$ model. We believe this is a new result not reported in the literature before.

Last, we presented fluid and diffusion approximations to a model of scheduled arrivals where the realized arrival epochs are uniformly random around the scheduled arrival epochs. We demonstrate that this model can be viewed as a special case of the $\D_{(i)}$ model.


While some non-stationary queueing models can be analyzed by the \textit{uniform acceleration} (UA) technique, we don't yet know if this can be used in analyzing `transitory' queueing models as well. Instead, we use the \textit{population acceleration} (PA) technique. An interesting question is the precise relationship between PA and UA techniques. UA is related to the notion of pointwise stationary approximation (PSA) of the performance metrics of the queue at a fixed time by stationary/ergodic states of an associated notional chain at that time. It is unclear to us if this can be done when pointwise ergodicity is not available. 


Finally, we consider this work to be a step towards a comprehensive `theory of transitory stochastic networks'. Our focus is next going to be to extend the current theory to `transitory queueing networks' with networks of $\D_{(i)}/GI/1$ queues. Moreover, obtaining stochastic process limits in the \textit{non-degenerate slowdown} regime \cite{At2012} (i.e., when the number of servers $s_n = n^{\a}$ for $\a \in [0,1]$) would be very useful in developing ``non-conventional'' large population approximations to transitory many-server systems. Another interesting class of problems is in the development of large and moderate deviation analyses of transitory queueing systems, something that we believe will be quite different from conventional queues. Last, while we have mostly focused on the single class setting, there are very interesting questions in the context of multi-class transitory systems, for instance in designing optimal schedulers.

\begin{APPENDICES} 
\section{Proofs of Lemmas and Theorems}
\subsection{Proof of Lemma \ref{lem:d-lim}}
\noindent (i). $\sD$ is a subspace under the relative topology $\t_r := \{A
\cap \sD | A \in \t_{\lim}\}$, where $\t_{\lim}$ is the topology
induced by the $J_1$ metric on $\sD_{\lim}$. Then, it follows that $\mathscr{D} =
\mathscr{D}_r := \{\sD \cap A | A \in \mathscr{D}_{\lim}\}$. To see
this, note that $\mathscr{D}_r$ contains all possible open sets of $\sD$
(since these are elements of $\mathscr{D}_{\lim}$) and $\mathscr{D}_r$
is a $\s$-algebra. This implies that $\mathscr{D} \subseteq
\mathscr{D}_r$ since the Borel $\s$-algebra is the smallest to contain
all open subsets of $\sD$. In the opposite direction, by since $\sD$
is a subspace, the injection map $\iota : \sD \to \sD_{\lim}$ is a
homeomorphism, implying that for any Borel set $A \in
\mathscr{D}_{\lim}$, $\iota^{-1}(A)$ is Borel in $\sD$. But, the
inclusion map is clearly $\iota^{-1}(A) = A \cap \sD$ by definition. Therefore
$\mathscr{D}_r \subseteq \mathscr D$.

This implies that for any $D \in \mathscr{D}_{\lim}$, $\mathbb P(x \in D) = \mathbb P(x \in D \cap \sD)$ is well defined since
$x$ is defined with respect to $(D,\mathscr{D})$. This
extends the definition of the measure induced by $x$ to
$\sD_{\lim}$. 

\noindent (ii) Now, let $G \subset \sD_{\lim}$ be any closed subset. Then, we have
\[
\bbP(x_n \in G) = \bbP(x_n \in G \cap \sD) + \bbP(x_n
\in G \cap \sD^c) = \bbP(x_n \in G \cap \sD),
\]
~where $\sD^c$ is the complement set. Let $E := G \cap \sD$, and note
that $\sD \backslash E = \sD \cap E^c = \sD \cap G^c$. This implies that $\sD
\backslash E$ is open in the relative topology on $\sD$ since $G^c$ is open in
$\t_{\lim}$. Now, using the fact
that $x_n \Rightarrow x$ in $(\sD,J_1)$ and part (iii) of
Theorem 2.1 of \cite{Bi68} we have $\limsup_{n \to \infty} \bbP(x_n \in G \cap \sD) \leq \bbP(x \in G \cap \sD)$. This implies that
\[
\limsup_{n \to \infty} \bbP(x_n \in G) \leq \bbP(x \in G \cap
\sD) = \bbP(x \in G),
\]
where the last equality follows by the fact that $x$ concentrates on
$\sD$. By part (i) of Theorem 2.1 of \cite{Bi68}, it follows that $x_n
\Rightarrow x$ in $(\mathcal D_{\lim},J_1)$ as $n \to \infty$.
\EndProof

\subsection{Statement and Proof of Lemma \ref{lem:X-hat}}
\begin{lemma}
\label{lem:X-hat}
~As $n \rightarrow \infty$,
\begin{equation}
\label{X-hat}
\hat{X}^n \Rightarrow \hat{X} \,:=\, \tilde W - W \circ \bar B ~ \text{ in } (\sD_{\lim},J_1)
\end{equation}
where $\bar{B}$ is defined in \eqref{busy-time-fluid}, and $\tilde W$ and
$W$ are mutually independent Gaussian bridge and Brownian motion
processes respectively, as defined in Proposition \ref{proposition:basics}
\end{lemma}

\Proof
Fix $j \in \{1, \ldots, N\}$. Recall that $B_j^n(t) \leq t, \forall t
\in [0,\infty)$,  implying that $S_j^n \circ B_j^n \in
\sD_{\lim}$. Using \eqref{lim:diffusion-ASV} in Proposition  \ref{proposition:basics}, Corollary \ref{cor:busy-time-fluid} and the random time change theorem (see, for example, Section 17 of \cite{Bi68}), it follows that
\(
\sqrt{N_n} \bigg ( \frac{S_j^n \circ B_j^n}{N_n} - \mu B_j^n\bigg ) \Rightarrow
\sigma \mu^{3/2} W_j \circ \bar{B}_j.
\)
~Now, it follows from Proposition \ref{proposition:basics} and the weak
limit just proved that
\(
\hat{X}^n \Rightarrow \hat{X}(t) \,:=\, \tilde W  - \sum_{j=1}^s \sigma
\mu^{3/2} W_j \circ \bar{B}_j \stackrel{d}{=} \tilde W - W \circ \bar B,
\)
where the final equality follows from the fact that the sum of independent
Brownian motions is equal in distribution to a Brownian motion.
\EndProof

\subsection{Statement and Proof of Lemma~\ref{lem:mama95}.}
This result is a consequence of Theorem 9.5.1 of
\cite{Wh01b}. A version of this result is also proved in \cite{HoJaWa12}.
\begin{lemma}[Directional derivative reflection mapping lemma] \label{lem:mama95}
~Let $x \in \sD$ and $y \in \sC$ be real-valued functions on $[0,\infty)$, and $\Psi(z)(t) = \sup_{0 \leq s \leq t} (-z(s))$, for any process $z \in \sD_{\lim}$. Let $\{y_n\} \subset \sD_{\lim}$ be a sequence of functions such that $y_n \stackrel{a.s.}{\rightarrow} y$ as $n \rightarrow \infty$. Then, with respect to Skorokhod's $M_1$ topology,
\(
\ty_n := \Psi(\sqrt{N_n} x + y_n) - \sqrt{N_n} \Psi(x) \longrightarrow \ty := \sup_{s \in \nabla_t^{x,L}}(-y(s)) \vee \sup_{s \in \nabla_t^{x,R}}(-y(s))
\)
as $n \rightarrow \infty$, where $\nabla_t^{x,\cdot}$ are defined in
Definition \ref{def:dirdev}.
\end{lemma}

  Rewrite $\ty_n$ as
\[
\ty_n = (\Psi(\sn x + y_n) - \Psi(\sn x + y)) - (\Psi(\sn x + y) - \sn \Psi(x)).
\]
Now, using the fact that the Skorokhod reflection map is Lipschitz continuous under the uniform metric (see Lemma 13.4.1 and Theorem 13.4.1 of \cite{Wh01}) we have
\[
(\Psi(\sn x + y_n) - \Psi(\sn x + y)) \leq \| y_n - y \|,
\]
where $\| \cdot \|$ is the uniform metric. It follows that
\[
\ty_n \leq \| y_n - y \| + (\Psi(\sn x + y) - \sn \Psi(x)),
\]
Now, by Theorem 9.5.1 of \cite{Wh01b} we know that as $n \rightarrow \infty$
\[
(\Psi(\sn x + y) - \sn \Psi(x)) \stackrel{a.s.}{\rightarrow} \ty, \text{ in } (\sD_{\lim},M_1).
\]
Using this result, and the fact that by hypothesis $y_n$ converges to $y$ in $(\sD_{\lim},U)$ we have
\[
\ty_n \stackrel{a.s.}{\rightarrow} \ty, \text{ in } (\sD_{\lim},M_1).
\]
\EndProof

\subsection{Proof of Proposition \ref{thm:j1-convergence}.} \label{sec:j1-convergence}
Recall that $\ty_n$ converges to $\ty$ in the $M_1$ topology, from Lemma \ref{lem:mama95}. Now, consider a path of $y$ that is non positive at $\t$. Thus, the
limit process $\ty$ has a discontinuity at $\t$ such that $\ty(\t) >
\ty(\t+)$. Note that the process path is left continuous at $\t$.
Assume that $\ty(\t) - \ty(\t+) > \d > 0$, and it follows that
$\ty(\t) = -y(\t) > \d$ (since $\ty(\t+) = 0$). Fix an $\e > 0$ such
that $\d > \e$. Now, by the continuity of $y$, there
exists $\eta > 0$ such that $\sup_{t \in [\t-\eta,\t+\eta]}
|y(t) - y(\t)| \leq \frac{\e}{4}$. Then, there also exists
a $n_0$ such that for all $n > n_0$, $0 \geq - \sn x(t) >
-\frac{e}{4}$ for $t \in [\t - \eta,\t]$. 

Then, for any $t \in [\t-\eta]$ it follows that
\[
-\sn x(t) - y(t) + y(\t) > - \frac{\e}{2}.
\]
This implies that
\(
-\sn x(t) - y(t) > \d - \frac{\e}{2} > \frac{e}{2},
\) since $\d > \e$. It follows that $\Psi(\sn x + y) >
\frac{\e}{2}$ for all time points $t \in [\t -\eta,\t]$. Thus, it
cannot be the case that uniform convergence is possible on any compact
set of $[-T_0,\infty)$. Furthermore, consider any sequence $\{\l_n\}
\subseteq \L$. Then, for large $n$, by assumption, $\l_n$ is uniformly
close to the identity map. Thus, any distortion introduced by the
homeomorphism will be minimal, and the same argument will show that it
cannot be the case that, for any fixed $\e > 0$, $\|\ty_n \circ \l_n -
\ty \circ e\| \leq \frac{\e}{2}$ for large $n$, and there is a set of points
determined by $\eta$ (due to the continuity of $y$) where it is
the case that $|(\ty_n \circ \l_n)(t) - (\ty \circ e) | > \frac{\e}{2}$.
\EndProof

\subsection{Proof of Lemma \ref{lem:avg-dist-limit}}
For each $n \in \bbN$ we
have $F_{n,i} = \Upsilon \left(  \frac{i}{N_n}\right)$, $i = 1,
\ldots, N_n$. Therefore, \eqref{eq:avg-dist} can be rewritten as
\(
\bar F_n (t) = \frac{1}{N_n}\sum_{i=1}^{N_n} \Upsilon \left(  \frac{i}{N_n}\right)(t),
\)
and we prove that
\(
\left \| \frac{1}{N_n}\sum_{i=1}^{N_n} \Upsilon \left(
  \frac{i}{N_n}\right)(t) - \int_{[0,1]} \Upsilon(s) (t) m(ds) \right\|_{[0,1]} \to 0
\)
 as $n \to \infty$. Notice that $\bar F(t)$ is a Riemann-Stieltjes integral with respect
to the Lebesgue measure. Therefore, it is natural to view $\bar F_n$
as a Riemann-Stieltjes (pre-limit) sum. For a fixed $t \in
[0,1]$, therefore, we show that the Riemann sums converge to the
Riemann-Stieltjes integral.

Let $M_i(t) := \sup \Upsilon(x)(t)$ and $m_i (t) :=
\inf \Upsilon(x)(t)$ for all $x \in \left[\frac{i-1}{N_n},
  \frac{i}{N_n}\right]$ and $i = 1, \ldots, N_n$. We define the ``upper'' and
``lower'' Riemann sums as (respectively) 
\(
U_n(t) := \sum_{i=1}^{N_n} M_i \left( \frac{i}{N_n} - \frac{i-1}{N_n}\right)
\)
 and
\(
N_n(t) := \sum_{i=1}^{N_n} m_i \left( \frac{i}{N_n} - \frac{i-1}{N_n}\right).
\)
Clearly, for every $\e > 0$ and large enough $n$, $U_n(t) - N_n(t) < \e$ due
to the Lipschitz continuity property assumed for $\Upsilon$. This is
tantamount to showing that the bound holds for at least one possible
partition of $[0,1]$. Then, by Theorem 6.6 of \cite{Ru1964}, it
follows that the limit exists and is equal to $\bar F(t)$. The
Lipschitz continuity property implies that the limit
clearly holds for all $t \in S$, implying uniform convergence.
\EndProof

\subsection{Proof of Lemma \ref{lem:avg-cov-limit}}
The first summation in the definition of $K_n(s,t) = \frac{1}{N_n} \sum_{i=1}^{N_n} F_{n,i}(s \wedge
t) - \frac{1}{N_n} \sum_{i=1}^{N_n} F_{n,i}(s) F_{n,i}(t)$ converges
to $\int_{\sK} F_p(s \wedge t) m(dp) $ as $n \to \infty$ by Lemma
\ref{lem:avg-dist-limit}. The second summation converges as well by
using the same Riemann-Stieltjes summation argument used in the proof
of Lemma \ref{lem:avg-dist-limit}, and the limit is $\int_{\sK} F_p(t)
F_p(s) m(dp)$. 
\EndProof

\subsection{Proof of Theorem \ref{thm:arrival-slln}}
First, fix $t \in [0,1]$ and $\e >0$. Consider
\begin{eqnarray}
\nonumber
\bbP \left(\left | \bar A^n(t) - \bar F_n(t) \right | > \e \right) &=&
\bbP \left(\left | \frac{1}{N_n}
  \sum_{i=1}^{N_n} \left( \mathbf 1_{\{T_{n,i} \leq t\}} -
    \Upsilon(i/N_n)(t)\right ) \right | > \e \right)\\
\nonumber
&\leq& \frac{1}{\e^4 N_n^4} \bbE \left | \sum_{i=1}^{N_n} \left( \mathbf
  1_{\{T_{n,i} \leq t\}} - \Upsilon(i/N_n)(t) \right) \right|^4\\
\nonumber
&=& \begin{split} &\frac{1}{\e^4 N_n^4} \sum_{i=1}^{N_n} \bbE \left| \mathbf
  1_{\{T_{n,i} \leq t\}} - \Upsilon(i/N_n)(t) \right|^4\\ &+ \frac{12}{\e^4 N_n^4} 
\sum_{i=1}^{N_n-1} \sum_{j=i+1}^{N_n} \bbE \left| \mathbf
  1_{\{T_{n,i} \leq t\}} - \Upsilon(i/N_n)(t)\right|^2 \bbE \left| \mathbf
  1_{\{T_{n,j} \leq t\}} - \Upsilon(j/N_n)(t)\right|^2  \end{split}\\
\label{eq:fourth-moment-bound}
&\leq& \begin{split}\frac{1}{\e^4 N_n^4} \sum_{i=1}^{N_n} \Upsilon(i/N_n)(t) &(1 -
\Upsilon(i/N_n)(t))\\ &+ \frac{12}{\e^4 N_n^4} \left( \sum_{i=1}^{N_n}  \Upsilon(i/N_n)(t) (1 -
\Upsilon(i/N_n)(t)) \right)^2,\end{split}
\end{eqnarray}
where the last inequality follows due to the fact that the terms that
remain in the expansion of $\bbE \left | \sum_{i=1}^{N_n} \left( \mathbf
  1_{\{T_{n,i} \leq t\}} - \Upsilon(i/N_n)(t) \right) \right|^4$ are
$\sum_{i=1}^{N_n} \bbE |\mathbf 1_{\{T_{n,i} \leq t\}} - \Upsilon(i/N_n)(t)|^2$ and
cross products $\sum_{i=1}^{N_n-1} \sum_{j=i+1}^{N_n} \bbE |\mathbf
1_{\{T_{n,i} \leq t\}} - \Upsilon(i/N_n)(t)|^2 \bbE |\mathbf
1_{\{T_{n,j} \leq t\}} - \Upsilon(j/N_n)(t)|^2$. From Lemma
\ref{lem:avg-cov-limit} it follows that \eqref{eq:fourth-moment-bound}
is bounded above by $\frac{C}{\e^4 N_n^2}$. Therefore, by the
Borel-Cantelli Lemma,
\(
\sum_{n=1}^{\infty} \bbP \left(\left | \bar A^n(t) - \bar F_n(t) \right
  | > \e \right) < \infty,
\)
implying that $\bbP \left(\left | \bar A^n(t) - \bar F_n(t) \right | >
  \e~ i.o. \right) = 0$. Combining this result with Lemma
\ref{lem:avg-dist-limit} proves that $\bar A^n$ converges to $\bar F$
almost surely pointwise. 

Next, consider a uniform partition of the support $[0,1]$, and suppose
$\frac{j-1}{M} \leq t \leq \frac{j}{M}$, where $j = 1, \ldots, M$ and
$M$ is the size of the partition. Then, for fixed $n$,
\(
\bar{A}^n\left( \frac{j-1}{M}\right) \leq \bar{A}^n\left(
  t\right) \leq \bar{A}^n\left( \frac{j}{M}\right), 
\)
implying that $\frac{1}{N_n} \sum_{i=1}^{N_n} \left( \mathbf 1_{\{T_{n,i} \leq j-1/M\}} -
F_{n,i}(j-1/M) \right)$
\[
\begin{split}
&\leq \frac{1}{N_n} \sum_{i=1}^{N_n} \left( \mathbf 1_{\{T_{n,i} \leq t\}} -
F_{n,i}(t) \right) + \frac{1}{N_n} \sum_{i=1}^{N_n} \left( F_{n,i}(j/M)
- F_{n,i}(j-1/M) \right)\\ &\leq \frac{1}{N_n} \sum_{i=1}^{N_n} \left( \mathbf 1_{\{T_{n,i} \leq j/M\}} -
F_{n,i}(j/M) \right).
\end{split}
\]
For each $M$, there exists $n_M$ such that for all $n \geq n_M$ $\left | F_{n,i}(j/M) -
  F_{n,i}(j-1/M) \right | \leq \frac{1}{M}$. Further, for $\e > 0$,
there exists $n^{'}_M$ such that for all $n \geq \max(n_M, n_M^{'})$,
\[
\left | \frac{1}{N_n} \sum_{i=1}^{N_n} \left( \mathbf 1_{\{T_{n,i} \leq k\}} -
F_{n,i}(k) \right) \right| < \e,
\]
where $k = j-1/M \text{ or } j/M$. It follows that
\[
\sup_{t \in [0,1]} \left | \frac{1}{N_n} \sum_{i=1}^{N_n} \left( \mathbf 1_{\{T_{n,i} \leq t\}} -
F_{n,i}(t) \right) \right | < 2(\e + \frac{1}{M}).
\]
Since $\e$ is arbitrary, letting $M \to \infty$ the desired result follows.
\EndProof

\subsection{Proof of Theorem \ref{thm:arrival-fclt}}
We first prove pointwise convergence by verifying the sufficiency of
the Lyapunov Central Limit Theorem (Theorem 7.3 \cite{Bi68}). Fix $t
\in [0,1]$ and let $\d > 0$, and consider
\[
\frac{\sum_{i=1}^{N_n} E|F_{n,i}(t) - \mathbf 1_{\{T_{n,i} \leq
    t\}}|^{2+\d}}{\left( \sum_{i=1}^{N_n} F_{n,i}(t) (1 - F_{n,i}(t)) \right)^{2+\d}}.
\]
Dividing the numerator and denominator by $1/N_n^{2+\d}$, note that
the denominator converges to $(\bar F(t))^{2+\d}$ as a consequence of
Lemma \ref{lem:avg-dist-limit}. Consider the numerator alone,
\[
\frac{1}{N_n^{2+\d}} \sum_{i=1}^{N_n} E|F_{n,i}(t) - \mathbf
1_{\{T_{n,i} \leq t\}}|^{2+\d} \leq \frac{2^\d}{N_n^{2+\d}}
\sum_{i=1}^{N_n} F_{n,i}(t) ( 1- F_{n,i}(t)),
\]
which tends to 0 as $n \to \infty$. The Lyapunov CLT implies that
\[
\frac{\sum_{i=1}^{N_n} (F_{n,i}(t) - \mathbf 1_{\{T_{n,i} \leq t\}})}{
  \sqrt{\sum_{i=1}^{N_n} F_{n,i}(t) (1 - F_{n,i}(t))}} =
\frac{\sum_{i=1}^{N_n} (F_{n,i}(t) - \mathbf 1_{\{T_{n,i} \leq
    t\}}}{\sqrt N_n} \times \frac{\sqrt{N_n}}{\sqrt{\sum_{i=1}^{N_n}
    F_{n,i}(t) (1 - F_{n,i}(t))}} \Rightarrow \sN(0,1).
\]
By Lemma \ref{lem:avg-cov-limit} it follows that $\frac{\sum_{i=1}^{N_n} (F_{n,i}(t) - \mathbf 1_{\{T_{n,i} \leq
    t\}}}{\sqrt N_n} \Rightarrow \tilde W(t) := \sqrt{K(t,t)} \sN(0,1)$. Next, using the Cramer-Wold device
it is straightforward to argue that $(\hat A^n(t_1), \ldots, \hat
A^n(t_k)) \Rightarrow (\tilde W(t_1), \ldots, \tilde W(t_k))$ where $(t_1, \ldots,
t_k) \in [0,1]^k$ for all $k \in \bbN$.

Finally, we verify the sufficiency of Theorem 15.6 of \cite{Bi68} to
show that $\hat A^n \Rightarrow \tilde W$ in $(\sD,J_1)$. To ease the notation, let $X_{n,i}(t) := (\mathbf 1_{\{T_{n,i} \leq
  t\}} - F_{n,i}(t))$. By Chebyshev's inequality, for any $\l > 0$ and $t_1 \leq
t \leq t_2 \in [0,1]$, $\l^4 \bbP(|\hat A^n(t) - \hat A^n(t_1)| \geq \l, |\hat A^n(t) - \hat
A^n(t_2)| \geq \l) $
\begin{eqnarray*}
&\leq& E \left[  (\hat A^n(t_1) - \hat A^n(t))^2
  (\hat A^n(t_1) - \hat A^n(t))^2\right].\\
&=& \frac{1}{N_n^2} E \left[ \left| \sum_{i=1}^{N_n} (X_{n,i}(t) -
X_{n,i}(t_1)) \right|^2 \left| \sum_{i=1}^{N_n} (X_{n,i}(t_2) -
X_{n,i}(t)) \right|^2 \right]\\
&=& \frac{1}{N_n^2} E \left[ \begin{split} &\sum_{i=1}^{N_n} |X_{n,i}(t) -
  X_{n,i}(t_1)|^2 \sum_{i=1}^{N_n} |X_{n,i}(t_2) - X_{n,i}(t)|^2\\
&+ 2 \sum_{i < j} (X_{n,i}(t_2) -X_{n,i}(t))(X_{n,j}(t_2) -
X_{n,j}(t)) \sum_{l=1}^{N_n} |X_{n,l}(t) - X_{n,l}(t_1)|^2\\
&+2 \sum_{i < j} (X_{n,i}(t) -X_{n,i}(t_1))(X_{n,j}(t) -
X_{n,j}(t_1)) \sum_{l=1}^{N_n} |X_{n,l}(t_2) - X_{n,l}(t)|^2\\
&+ 4 \sum_{i < j} (X_{n,i}(t) - X_{n,i}(t_1))(X_{n,j}(t) -
X_{n,j}(t_1)) \sum_{i < j} (X_{n,i}(t_2) - X_{n,i}(t))(X_{n,j}(t_2) -
X_{n,j}(t)) \end{split} \right] \\
&\leq& 
\begin{split}
&\frac{1}{N_n} \sum_{i=1}^{N_n} \left[ (F_{n,i}(t) -
  F_{n,i}(t_1))(F_{n,i}(t_2) - F_{n,i}(t))(1 - F_{n,i}(t) +
  F_{n,i}(t_1))(1 - F_{n,i}(t_2) + F_{n,i}(t))\right]^{1/2}\\
& + 2 \frac{1}{N_n^2} \sum_{i<j}  \left[ (F_{n,i}(t) -
  F_{n,i}(t_1))(F_{n,j}(t_2) - F_{n,j}(t))(1 - F_{n,i}(t) +
  F_{n,i}(t_1))(1 - F_{n,j}(t_2) + F_{n,j}(t))\right]\\
&+ 4 \sum_{i < j} (F_{n,i}(t) - F_{n,i}(t_1))(F_{n,i}(t_2) -
F_{n,i}(t))(F_{n,j}(t) - F_{n,j}(t_1))(F_{n,j}(t_2) - F_{n,j}(t))
\end{split}\\
&\leq& C
\end{eqnarray*}
where $C \geq 8$ and the bound is true for all $t_2 >
t_1$. Theorem 15.6 of \cite{Bi68} shows that if 
\(
\bbP(|\hat A^n(t) - \hat A^n(t_1)| \geq \l, |\hat A^n(t) - \hat
A^n(t_2)| \geq \l) \leq (G(t_2) - G(t_1))^{2\a},
\)
where $G$ is a non-decreasing function on $[0,1]$ and $\a > 1/2$, then
$\hat A^n$ converges weakly to a limit in $(\sD,J_1)$. Therefore, $\hat
A^n \Rightarrow \tilde W$ as $n \to \infty$. The convergence in
$(\sD_{\lim},J_1)$ follows by an application of part (ii) of Lemma
\ref{lem:d-lim}. Finally, by part (ii) of Theorem 1.1 of \cite{Sh1979} we know that
$\tilde W$ has continuous sample paths, implying that $A^n \Rightarrow
\tilde W$ in $(\sD_{\lim}, U)$, thus completing the proof.
\EndProof

\subsection{Proof of Theorem \ref{thm:conditioned-poisson-fclt}}
Let $\mathbf T := \{T_i, \, i=1,\ldots,n\}$ be a collection of i.i.d. random
variables, with distribution function $F$ (defined in
\eqref{eq:dist-func-from-poisson}). Let $A^n(t) := \sum_{i=1}^n
\mathbf{1}_{\{T_i \leq t\}}$ and $\hat A^n(t) := \sn
\left( \frac{A^n(t)}{n} - F(t) \right)$ be the
empirical process associated with $\mathbf T$. 
For a fixed $n \geq 1$ and $x \in \bbR$ we have
\[
\bbP(\hat M_n (t) \leq x | M(T) = n) = \bbP(M(t) \leq x \sn + n F(t) |
M(T) = n).
\]
The \textsc{OS} property implies that $M(t) |_{\{M(T) = n\}}
\stackrel{d}{=} A^n(t)$. Proposition \ref{prop:empirical-process-limit}
implies that 
\begin{eqnarray*}
\bbP(M(t) \leq x \sn + n F(t) | M(T) = n) &=& \bbP(A^n(t) \leq x \sn + n
F(t))\\
&=& \bbP(\hat A^n(t) \leq x) \\
&\Rightarrow& (W^0 \circ F)(t)
\end{eqnarray*}
proving the pointwise convergence of the process $\hat M_n$. It is
also well known that for any $0 < t_1 < \cdots < t_d < T$,
$\bbP(M(t_1) = n_1, \ldots, M(t_d) = n_d | M(T) = n) = \bbP(A^n(t_1) =
n_1, \ldots, A^n(t_d) = n_d)$, where $n_1 + \cdots n_d = n$, so the fact
that $(\hat A^n(t_1), \ldots, \hat A^n(t_d))
\Rightarrow (W^0 \circ F)(t_1), \cdots, (W^0 \circ F)(t_d))$ implies
that the finite dimensional distrbutions of $\hat M_n$ converge to the
same limit. The tightness of $\hat M_n$ is implied directly by that of $\hat A^n$,
so by Theorem 8.1 of \cite{Bi68} the theorem is proved.
\EndProof

\subsection{Proof of Proposition \ref{prop:exchangeinter}}
Let $\{x_1, \ldots, x_n\} \subset [0,T]$ be such that $0 \leq x_1 < x_2 < \cdots < x_n \leq T$. Consider the measure of the event $\{\xi_{1} \in dx_1, \cdots, \xi_{n} \in dx_n\} \in \sF$,
\[
\begin{split}
P(\xi_{1} \in dx_1,& \cdots, \xi_{n} \in dx_n | M(T) = n) \\= 
&\frac{P((\xi_{1} \in dx_1, \cdots, \xi_{n} \in dx_n), M(T) =n)}{P(M(T) = n)}.
\end{split}
\]
Recall that $\{M(T) = n\} = \{\sum_{l=1}^n \xi_{l} \leq T < \sum_{l=1}^{n+1} \xi_{l}\}$, implying that we have:
\[
\begin{split}
&\frac{P((\xi_{1} \in dx_1,\cdots, \xi_{n} \in dx_n), M(T) =
  n)}{P(M(T) = n)} \\= 
&\frac{ P(\xi_{1} \in dx_1, \ldots \xi_{n} \in dx_n,\sum_{l=1}^n
  \xi_{l} \leq T, \sum_{l=1}^n\xi_{l} + \xi_{n+1} > T)}{P(M(T) =n)}. 
\end{split}
\]
Using the fact that under the measure $\bbP$, $\xi_{i}$ are i.i.d. random variables, it follows that the measure of the joint event is invariant under any permutation of the first $n$ random variables. That is, if $\pi(\cdot)$ is a permutation of $\{1,\ldots,n\}$, then we have
\[
\begin{split} 
& \frac{P((\xi_{1} \in dx_1, \cdots, \xi_{n} \in dx_n), M(T) =
  n)}{P(M(T) = n)}\\ = &\frac{P(\xi_{\pi(1)} \in dx_1, \ldots
  \xi_{\pi(n)} \in dx_n,\sum_{l=1}^n \xi_{\pi(l)} \leq T,
  \sum_{l=1}^n\xi_{\pi(l)} + \xi_{n+1} > T) }{P(M(T) = n)},
\end{split}
\]
which is equal to $P_n(\xi_{\pi(1)} \in dx_1, \ldots \xi_{\pi(n)} \in dx_n)$. Next, suppose that  $\pi(\cdot)$ is a permutation of $\{1,\ldots, n+1\}$. Then, it is possible that the event $\sum_{i=1}^n \xi_{\pi(l)} > T$, since $\xi_{n+1} > T - \sum_{l=1}^n \xi_{l} > 0$, conditionally on $\{M(T) = n\}$. Thus, $\Xi_n$ cannot be extended to a larger collection of exchangeable random variables, implying that it is finitely exchangeable.
\EndProof

\subsection{Conditioned Renewal Process: Lemmata} \label{sec:lemmata}
\subsubsection{Sample Space Construction}
We assume that the underlying sample space $\Omega, \sF, \bbP$ is rich
enough to support a sequence of (jointly) independent stochastic
processes $\{M_n\}, \, n \geq 1$, such that they are each
\textit{indistinguishable} from $M$. That is, for any $n \geq 1$, $\bbP(M_n(t)
\not = M(t),\, \forall t \geq 0) = 0$. For a fixed $n \geq 1$ and $T
> 0$, we define the restricted sample space, $(\Omega_n, \sF_n,
\bbP_n)$, where $\Omega_n = \Omega \cap \{M_n(T) = n\}$, $\sF_n := \sigma \{A
\cap \{M_n(T) = n\} : A \in \sF\}$ and $\bbP_n(B) :=
\frac{P(B)}{P(M_n(T)= n)}$ for any $B \in \sF_n$. Clearly $\{ \{M(T) =
n\}$ $n \geq 1 \}$ forms a partition of $\Omega$. The following claim shows that this property ($\bbP-a.s.$)
extends to the collection $\{\Omega_n\}$.

\begin{lemma} \label{lem:partition}
~$\{\Omega_n\}$ $n \geq 1$ forms a $\bbP-a.s.$ partition of $\Omega$.
\end{lemma} 
\Proof
For a fixed $n \geq 1$ and $m \in \bbN$, $\{M_n(T) = m\}$, forms
a partition of $\Omega$, as do $\{M(T) = m\}$. By assumption $M_n$ and $M$ are
indistinguishable from each other. It is straightforward to deduce
that $\bbP(\{M_n(T) = m \} \Delta \{M(T) = m\}) = 0$, implying that
$\{M_n(T) = m\} = \{M(T) = m\} \, \bbP-a.s.$ for every $m \in \bbN$. 

Now, consider the collection of sets $\{ \{M_n(T) =
 n\} \}$ $n \geq 1$. For brevity, let $A_n := \{M_n(T) = n\}$ and $B_n
 := \{M(T) = n\}$. We have
\begin{eqnarray*}
\left(\cup_{n \geq 1} \{M_n(T) = n\} \right) \Delta \left(\cup_{n \geq
  1} \{M(T) = n\}\right)  &=& \left( \cup_{n \geq 1} A_n \right) \Delta
\left( \cup_{n \geq 1} B_n\right) \\
&=& \cup_{l \geq 1} \left( \cap_{n \geq 1} (B_l \cap A_n^c)\right).
\end{eqnarray*}
By the assumption of indistinguishability, $\bbP(B_l \cap A_l^c) =
0$, implying that $\bbP(\cap_{n \geq 1} (B_l \cap A_n^c)) = 0$, for
every $l \geq 1$. Therefore, 
\[
\bbP \left(\cup_{n \geq 1} \{M_n(T) = n\} \right) \Delta \left(\cup_{n \geq
  1} \{M(T) = n\}\right) = 0.
\]
By virtue of the fact that $\cup_{n \geq 1} \{M(T) = n\} = \Omega$ it follows that $\{\{M_n(T) = n\} \}, \, n \geq 1$
forms a partition of $\Omega$ $\bbP-a.s.$
\EndProof

Next, we construct a new product space from the restricted sample
spaces $(\Omega_n, \sF_n, \bbP_n)$ as follows. Let, $\bar{\Omega} :=
\Omega_1 \times \Omega_2 \times \cdots$, so that $A \subset
\bar{\Omega} = A_1 \times A_2 \times \cdots$ for sets $A_n \subset
\Omega_n$. The product $\sigma-$algebra, $\bar{\sF} := \sF_1 \otimes
\sF_2 \otimes \cdots$ is the $\sigma-$algebra generated from cylinder sets
of the type $R = \{(\omega_1, \omega_2, \cdots) \in \bar\Omega |
\omega_{i_1} \in A_{i_1}, \cdots, \omega_{i_k} \in A_{i_k}\}$, where
$(i_1, \ldots, i_l)$ is an arbitrary subset of $\bbN$ of size $k \geq
1$ and $A_{i_n} \in
\sF_n$. The existence of such a product $\sigma-$algebra is well-justified by Proposition
1.3 in \cite{Fo1984}. Finally, we define $\bar{\bbP}(R) = \Pi_{i=1}^k
\bbP_{i_l}(A_{i_l})$, for the cylinder sets. This extends to $\bar{P}
= \bbP_1 \times \bbP_2 \times \cdots$, which is the natural product measure on the
measure space $(\bar{\Omega}, \bar{\sF})$, by standard arguments 
showing that the measure is countably additive on $\bar\sF$. The definition of the
Lebesgue integral on the space $(\bar{\Omega}, \bar{\sF}, \bar\bbP)$ now
follows from standard definitions of integration on product
spaces. However, we introduce some notation to ease our burden. In
particular, consider a function defined in the following manner: $\bar{X} := X \times \Pi_{l \not = n}
\bbI_{\{\Omega_l\}}$, where $X$ is measurable and
integrable with respect to $(\Omega_n, \sF_n, \bbP_n)$, and
$\bbI_{\{\cdot\}}$ is the indicator function. Then
\[
E_{\bar{\bbP}} [\bar{X}] = \int_{\Omega_n} X d\bbP_n \int_{\Pi_{l \not
  = n}} \bbI_{\{\Omega_l\}} d \bbP_l
\]
is well-defined, and we write this as $E_{\bar\bbP}[X]$, where it is
to be understood that the integration is actually of $\bar X$.

\subsubsection{Preliminary Lemmata}
Armed with the new product space, we can now proceed to the proof of
the diffusion limit. We define the
collection of random variables, $\Xi_{n} := (\xi_{n,1}, \cdots,
\xi_{n,n})$ to be ``conditioned'' inter-arrival times of the renewal process $M_n$
provided that
\[
M_n(t) = \sup \{k \geq 0 | \sum_{i=1}^k \xi_{n,i} \leq t\}, 
\]
for $t \leq T$ and $\sum_{i=1}^n \xi_{n,i} \leq T$.
From Proposition \ref{prop:exchangeinter} it follows that $\Xi_n$ forms
an exchangeable array, under the measure $\bbP_n$. Furthermore, under the
measure $\bar \bbP$, $\Xi_n$ and $\Xi_m$ ($m \not = n$) are
independent of each other (by the definition of the product
measure). The following lemma characterizes the mean and variance of
the inter-arrival times, for a fixed $n$. Intuitively, one should
expect the mean and variance to decrease to 0 as $n \to \infty$, as
there are a larger number of variables being packed into a fixed
interval ($T$ is fixed, of course). This characterization is important
in proving the diffusion limit. We start with a simple lemma.

\begin{lemma} \label{lem:integrable}
~Let $\{f_n\}$ be a sequence of non-negative, measurable,
functions. Let $P$ be a given measure. Assume: a) $\int f_n dP \to 0$
as $n \to \infty$, b) $|f_n| \leq C < \infty$. Then, $\lim_{n \to \infty} f_n = 0$ as $n \to \infty$.
\end{lemma}

\Proof
\begin{eqnarray*}
\lim \int f_n dP &\stackrel{DCT}{=}& \int \lim f_n dP\\ 
0 &=& \int \lim f_n dP,
\end{eqnarray*}
from whence the conclusion follows easily.
\EndProof

The proof of the asmptotic negligibility of the mean and variance
follow as consequences of Lemma \ref{lem:integrable}.

\begin{lemma} \label{lem:condition}
~For $\xi_{n,j} \in \Xi_n$,
\begin{eqnarray*}
(i) & \m_n := E[\xi_{n,j} | M_n(T) = n] \to 0 \text{ as } n \to \infty.\\
(ii) & E[|\xi_{n,j} - \m_n|^2 | M_n(T) = n] \to 0 \text{ as } n \to \infty.
\end{eqnarray*}
\end{lemma}

\Proof
\noindent (i) Conditioned on $\{M_n(T) = n\}$, $\{\xi_{n,j}\}_{j \leq n}$ are exchangeable, implying that they have the same distribution. Thus, 
\[
E[\xi_{n,1} | M_n(T) =n ] = \frac{1}{n} E[S_n | M_n(T) = n],
\]
where $S_{n} = \sum_{i=1}^n \xi_{n,i}$. By the definition of conditioned expectation:
\begin{eqnarray*}
\int \mathbf{1}_{\{M_n(T) = n\}} S_n d\bbP &=& \int \mathbf{1}_{\{S_n \leq T \leq S_{n+1}\}} S_n d\bbP  \\
&\leq& T \bbP(S_{n} \leq T < S_{n+1}) \leq T \bbP(S_n \leq T).
\end{eqnarray*}
Note: we interpret $S_{n+1}$ as the sum of $n+1$ inter-event times in
the $n$th unconditioned system. By assumption, $E_{\bbP}[\xi_{n,j}] = \mu > 0$, implying (by SLLN or Second Borel Cantelli Lemma) that $S_n \to \infty \, a.s.$ as $n \to \infty$. This implies that $\lim_{n \to \infty} \bbP(S_n \leq T) = 0$. Now, since $S_n > 0$, it follows from Lemma \ref{lem:integrable} that 
\[
\lim_{n \to \infty} E[S_n | M_n(T) =n] = 0.
\]
Finally, $\frac{1}{n} E[S_n|M_n(T) = n] \leq E[S_n | M_n(T) = n] \to 0$.\\

\noindent (ii) follows by a similar argument. From part (i) we have
convergence in measure, using Markov's inequality: for any $\e > 0$,
\begin{eqnarray*}
\bbP_n(\xi_n > \e) = \bbP(\xi_n > \e | M_n(T) = n) &\leq& \frac{E[\xi_n | M_n(T) = n]}{\e} \\
& \to & 0,
\end{eqnarray*}
as $n \to \infty$. Next, to see that the
variance converges to 0 as well, consider the following: recall that
$\bbP_n$ is the conditional probability measure (this is regular, in fact). Then, we have,
\begin{eqnarray*}
E[(\xi_{n,j} - \m_n)^2 | M_n(T) = n] &=& \int (\xi_{n,j} - \m_n)^2 d\bbP_n\\
&=& \int_{\{ \xi_{n,j} > \e\}} (\xi_{n,j} - \m_n)^2 d\bbP_n + \int_{\{\xi_{n,j} \leq \e\}} (\xi_{n,j} - \m_n)^2 d\bbP_n\\
& \leq &2 T^2 \bbP_n(\xi_{nj} > \e) + \e^2 \bbP_n(\xi_{nj} \leq \e) + \m_n^2 \\
& < & (2 T^2 + 2) \e
\end{eqnarray*}
for large $n$. Note that here, we used the fact that under the measure $\bbP_n$ the random variables are bounded by $T$, as well as the convergence in measure noted above. The same argument, in fact, extends to any $r$-th mean.
\EndProof
It is straightforward to see that Lemma \ref{lem:condition} is true
under measure $\bar \bbP$. The next theorem proves a rate of
convergence for the mean inter-arrival times in the conditioned
systems.

\begin{lemma} \label{thm:condition-mean-rate}
~Let $M$ be a renewal process with renewal time distribution $F$. Let $\m_n$ be the mean inter-arrival time, when the process is conditioned to have $n$ events by time $T$. Then, (i) $\m_n \to 0$ and (ii) $\sn |1 - n \m_n| \to 0$ as $n \to \infty$.
\end{lemma}
\Proof
Let $F: \bbR_+ \to [0,1]$ be the inter-arrival time distribution, defining a renewal process $M(t)$. Assume that $f(t) := \frac{d F(t)}{dt}$ is well defined. Then, the \emph{conditional intensity function} (CIF) of $M(t)$ is $\l^*(t) := E[N(dt)| \mathcal{H}_t] = \frac{f(t) dt}{1 - F(t)} \geq 0$, where $\mathcal{H}_t$ is the filtration generated by $M$. The integrated CIF, $\L^*(t) = \int_0^t \l^*(s) ds$ is the martingale \emph{compensator} so that $M(t) - \L^*(t)$ is a Martingale process. Theorem 7.4.I of \cite{DaVe2003a} shows that $M(t) = N(\L^{*^{-1}}(t))$ is a unit rate Poisson Process. That is, if $\{T_i\}$ is a realization of the event times of process $M$, then $\{\tT_i = \L^*(T_i)\}$ is a realization of those of a unit rate Poisson process. As $\L^*$ is non-decreasing, it follows that $\{\tT_{n+1} > \L^*(T) \geq \tT_n\}$ if and only if $\{T_{n+1} > T \geq T_{n}\}$, implying that $\{M(T) = n\} = \{N(\L^*(T)) = n\}$.

Let $\xi_i = T_i - T_{i-1}$ be the inter-arrival time random variable. Then,
\[
\bbP(\xi_1 > u | M(T) = n) = P(\phi_1 \geq \L^*(u) | N(\L^*(T)) = n),
\]
where $\phi_1 = \L^*(\xi_1)$. Consider the latter conditioned probability, and recall that a Poisson process satisfies the \textsc{Ordered Statistics} property. It follows that
\[
P(\phi_1 \geq \L^*(u) | N(\L^*(T)) = n) = \left( 1 - \frac{\L^*(u)}{\L^*(T)} \right)^n.
\]
As the inter-arrival times are exchangeable (when conditioned by $\{N(\L^*(T)) = n\}$) they are also identically distributed, so the ensuing arguments hold true for any inter-arrival time $\phi_i$ $i = 1, \ldots, n$. The conditional distribution implies that
\[
\tilde{\m}_n := \int_0^{\L^*(T)} P(\phi_1 \geq \L^*(u) | N(\L^*(T)) = n) d \L^*(u) = \frac{1}{n+1}.
\]
Equivalently, after a time change, we have
\[
\int_0^T \bbP(\xi_1 > u | M(T) = n) \l^*(u) du = \frac{1}{n+1}.
\]

However, we are interested in the asymptotics of the closely related integral $\m_n := \int_0^T \bbP(\xi_1 > u | M(T) = n) du$. Therefore, consider
\begin{eqnarray}
(n+1) \left| \m_n - \frac{1}{(n+1)} \right| &=& (n+1) \left| \int_0^T \bbP(\xi_1 > u | M(T) = n)(1 - \l^*(u)) du \right| \\
& \leq & K (n+1) \int_0^T \bbP(\xi_1 > u | M(T) = n) du\\
&=& K \int_0^T (n+1) \left( 1 - \frac{\L^*(u)}{\L^*(T)}\right)^n du,
\end{eqnarray}
where the inequality follows as the CIF is bounded on compact intervals. Since $0 \leq \frac{\L^*(u)}{\L^*(T)} \leq 1$ for every $0 < u \leq T$, it follows that $(n+1) \left( 1 - \frac{\L^*(u)}{\L^*(T)}\right)^n \to 0$ as $n \to \infty$. Then, using (the reverse) Fatou's Lemma we have
\[
\limsup_{n \to \infty} (n+1) \left| \m_n - \frac{1}{(n+1)} \right|  \leq \int_0^T \limsup_{n \to \infty} (n+1) \left( 1 - \frac{\L^*(u)}{\L^*(T)}\right)^n du = 0.
\]
Thus, Lebesgue almost everywhere on $[0,T]$, we have $\m_n \sim \frac{1}{n+1}$, so that $\m_n \to 0$ as $n \to \infty$. This immediately implies that $\sn |1 - n \m_n|\to 0$.
\EndProof

\subsection{Proof of Theorem \ref{thm:conditioned-fslln}}
Consider the interval $[0,1]$, and consider
\[
| \sum_{l=1}^{\lfloor nt \rfloor} \xi_{n,l} - t | \leq
|\sum_{l=1}^{\lfloor nt \rfloor} (\xi_{n,l} - \m_n) | + |\lfloor nt \rfloor \mu_n - t|.
\]
The second term on the RHS tends to $0$, as a consequence of Lemma
\ref{thm:condition-mean-rate} (it is straightforward to see that
$\{\xi_{n,l}\}_{l=1}^n$ satisfies the conditions of the theorem under
the space $(\bar \Omega, \bar \sF, \bar \bbP )$). For the first term,
consider the martingale sequence
\(
z_{n,l} := (\xi_{n,l} - \m_n) - E[(\xi_{n,l} - \m_n) | \sF_{n,l-1}],
\)
 where $\sF_{n,l} := \s \{(\xi_{n,1} - \m_n), \ldots, (\xi_{n,l-1} -
 \m_n), \sum_{i=l}^n (\xi_{n,i} - \m_n)\}$. Using the fact that
 $\{\xi_{n,l}\}$ are exchangeable it is easy to deduce that
\[
\sum_{i=j}^n (\xi_{n,i} - \m_n) = \sum_{i=j}^n E[(\xi_{n,i} - \m_n) |
\sF_{n,j-1}] = (n-j+1) E[\xi_{n,j} - \m_n | \sF_{n,j-1}]. 
\]
This implies that
\(
z_{n,l} = (\xi_{n,l} - \m_n) + \frac{1}{n-j+1} \sum_{i=l}^n (\xi_{n,i}
- \m_n) = \xi_{n,l} - \frac{1}{n-l+1} \sum_{i=l}^n \xi_{n,l}.
\)
 Using the fact that $\xi_{n,l} \in [0,1]$, under the measure $\bar
 \bbP$, it follows that
\begin{eqnarray*}
\sum_{l=1}^n \frac{1}{n-l+1} \sum_{i=l}^n \xi_{n,l} &\leq&
\sum_{l=0}^{n-1} \frac{1}{n-l} \sum_{j=l+1}^n \xi_{n,j}\\
&\leq& 2 \left( n - \sum_{l=0}^{n-1} \frac{1}{n-l} \right)\\
&\approx& 2 \left( n - \log n +o(n) \right).
\end{eqnarray*}

By definition we have, for any $\e > 0$,
\[
\{\omega \in \bar \Omega : |\sum_{l=1}^n \xi_{n,l} - \m_n|\} \subseteq
\{\omega \in \bar \Omega : |\sum_{l=1}^n z_{n,l}| > \e - n \m_n -2 (n
+ \log n +o(n)) \}.
\]
Thus, it suffices to bound the latter event, which we do using the
Azuma-Hoeffding inequality for Martingale differences.  First, let us recall the statement of the
Azuma-Hoeffding inequality. Let $Z_{l}$, $1 \leq l \leq n$, be a
martingale difference sequence defined with respect to a sample space $(\bbS, \sS, \sP)$, such that $|Z_l| \leq c_l$, for some set of constants $\{c_1, \ldots, c_n\}$. Then, for all $\e > 0$, 
\[
\sP\left(\max_{1\leq k \leq n} \sum_{l=1}^k Z_l > \e \right) \leq \exp \left( - \frac{\e^2}{2 \sum_{l=1}^n c_l^2}\right).
\]

Returning to the problem at hand, we have, for any $\e > 0$
\[
\bar \bbP (|\sum_{l=1}^n Z_{n,l}| > \e - n \m_n -2 (n
+ \log n +o(n))) \leq 2 \exp(- \frac{(\e - n \m_n -2 (n
+ \log n )^2}{4n}),
\]
where $|z_{n,l}| \leq 2$. Now, all that is required is a lower bound
on the exponentiated expression. Let
\[
f(\e,n) := \frac{1}{4n} \left( (\e - n \m_n -2 (n
+ \log n )^2\right).
\]
Multiplying and dividing by $n^2$ on the RHS, we obtain,
\begin{eqnarray*}
f(\e,n) &=& \frac{n}{4} \left( \e^2 +(\m_n + 2) (1 - 4 \e/n) + 8
  (\frac{\log n}{n})^2 + 4 \frac{\log n}{n} - 4 \e \frac{\log n}{n} (1
  - \e/n) \right)\\
&\geq& \frac{n}{4} \left( \e^2 + (1 - 4 \e/n)(1/n + 2 + 4 \frac{\log
    n}{n}) + 8 (\frac{\log n }{n})^2 \right),
\end{eqnarray*}
where in the last step we make use of the fact that $\frac{1}{n} >
\m_n $. For large enough $n$ such that $n > 4 \e$ we have
\[
f(\e,n) \geq \frac{n\e^2}{4} + \frac{2 \kappa_\e}{n},
\]
where $\kappa_{\e} := \log 4\e$. It follows that
\[
\exp (- f(\e,n)) \leq \exp(- \frac{n \e^2}{4})
\exp(-\frac{2\kappa_{\e}}{n}) \leq \exp(- \frac{n \e^2}{4}).
\]
This implies that for any $\e > 0$
\[
\sum_{n=1}^{\infty} \bar \bbP(|\sum_{l=1}^n(\xi_{n,l} - \m_n)| > \e) \leq
\sum_{n=1}^{\infty} \bar \bbP (|\sum_{l=1}^n Z_{n,l}| > \e - n \m_n -2 (n
+ \log n +o(n))) < \infty.
\]
Thus, by the First Borel-Cantelli Lemma, $\bar
\bbP(|\sum_{l=1}^n(\xi_{n,l} - \m_n)| > \e) i.o. ) = 0$. Therefore,
$\sum_{l=1}^n(\xi_{n,l} - \m_n) \to 0 $ a.s. as $n \to
\infty$. Clearly, this holds true for any $t \in [0,1]$, so that
\(
\sum_{l=1}^{\lfloor nt \rfloor}(\xi_{n,l} - \m_n) \to 0
\)
 a.s. as $n \to \infty$. By standard arguments it follows that the
 limit holds uniformly on $[0,1]$ concluding the theorem.
\EndProof

\subsection{Statement and Proof of Proposition \ref{prop:verification}}
\begin{proposition} \label{prop:verification}
~The triangular array $\{\phi_{n,l}, ~~ l=1,\ldots,n\}~~ n \geq 1$
satisfies the following properties:
\noindent (i) $\sum_{l=1}^n \phi_{n,l} \stackrel{P}{\to} 0$ as $n \to \infty$.\\
\noindent (ii) $\max_{1 \leq l \leq n} |\phi_{n,l}| \to 0$ as $n \to \infty$.\\
\noindent (iii) $\sum_{l=1}^n \phi_{n,l}^2 \stackrel{P}{\to} 1$ as $n
\to \infty$.
\end{proposition}
\Proof
~(i) The proof follows by using the definition of $\bar{\bbP}$. Fix $\e > 0$, and consider the following:
\begin{eqnarray*}
\bar{\bbP}(|\sum_{l=1}^n \phi_{n,l}| > \e) &=& \bar{\bbP}(|\sum_{l=1}^n \phi_{n,l}| > \e, M_n(T) = n)\\
&=& \bar{\bbP}(|\sum_{l=1}^n \xi_{n,l} - n \mu_n| > \e \sqrt{n}, M_n(T) = n)\\
&=& \begin{split} \bar{\bbP}&(\sum_{l=1}^n \xi_{n,l} > \e \sqrt{n} + n
  \mu_n , M_n(T) = n)\\  &+ \bar{\bbP}(\sum_{l=1}^n \xi_{n,l} < -\e
  \sqrt{n} + n \mu_n , M_n(T) = n). \end{split} 
\end{eqnarray*}
The first equality follows by the fact that under $\bar{\bbP}$ $\{M_n(
t) = n\}$, for every $n \geq 1$, are full measure sets. Recalling that
$\{M_n(T) = n\} = \{\sum_{l=1}^n \xi_{n,l} \leq T < \sum_{l=1}^n \xi_{n,l} + \xi_{n,n+1}\}$, it follows that for any $\omega \in A_n := \{\sum_{l=1}^n \xi_{n,l} > \e \sqrt{n} + n \mu_n , M_n(T) = n)\}$ we have
\[
T \geq \sum_{l=1}^n \xi_{n,l} > \e \sqrt{n} + n \mu_n.
\]
Now, using the fact that $\xi_{n,l}$ are exchangeable (for a fixed $n$), it follows directly that $n \mu_n = E_{\bar{\bbP}}[\sum_{l=1}^n \xi_{n,l}] \leq T$, under the measure $\bar{\bbP}$. Therefore, $n \mu_n$ is uniformly bounded (for every $n \geq 1$). It follows that for a given $T$, there exists a $n_T$ such that for every $n > n_T$, $\sqrt{n} \e + n \mu \geq T$. As $\e > 0$ is arbitrary, asymptotically, $A_n$ is an impossible event. Next, consider $B_n := \{\sum_{l=1}^n \xi_{n,l} < -\e \sqrt{n} + n \mu_n , M_n(T) = n\}$. Similar arguments show that 
\[
-\e \sqrt{n} + n \m_n > \sum_{l=1}^n \xi_{n,l} \geq 0.
\]
Clearly as $n \to \infty$, $-\e \sqrt{n} + n \m_n \to -\infty$, as $n \m_n$ is uniformly bounded. Since $\e > 0$ is arbitrary, $B_n$ too is (asymptotically) an impossible event. It follows that $\phi_{n,l} \stackrel{P}{\to} 0$ as $n \to \infty$.

\noindent (ii) The proof is elementary. First, using the union bound we have, for a fixed $\e > 0$,
\[
\bar{\bbP} (\max_{1 \leq l \leq n} |\phi_{n,l}| > \e) \leq \sum_{l=1}^n \bar{\bbP}(|\phi_{n,l}| > \e).
\]
Using the fact that the random variables $\phi_{n,l}$, $l \leq n$ are exchangeable, they are also (marginally) identically distributed. Thus,
\begin{eqnarray*}
\bar{\bbP} (\max_{1 \leq l \leq n} |\phi_{n,l}| > \e) &\leq& n \bar{\bbP}(|\phi_{n,1}| > \e)\\
&\leq& n \frac{E_{\bar{\bbP}} |\xi_{n,l} - \m_n|^2}{n\e^2 } = \frac{\sigma_n^2}{\e^2},
\end{eqnarray*}
where the latter expression follows by an application of Chebyshev's inequality under the $\bar{\bbP}$ measure. As noted before, $\sigma^2_n \to 0$ as $n \to \infty$. As $\e > 0$ is arbitrary, (ii) is proved.

\noindent (iii) The proof is more involved, requiring the construction
of a martingale difference sequence, and then an appeal to the
Azuma-Hoeffding inequality for a tight bound on the martingale
difference sequence. Note that $\phi_{n,l}^2$, $1 \leq l \leq n$, is an exchangeable sequence. Let $Z_{n,l} := \phi_{n,l}^2 - E_{\bar{\bbP}}[\phi_{n,l}^2 | \sF_{n,l-1}]$, where $\{\sF_{n,l}\}$ is a filtration defined with respect to $\phi_{n,l}^2$ as 
\[
\sF_{n,l} = \sigma( \phi_{n,1}^2, \ldots, \phi_{n,l-1}^2, \sum_{i=l}^n \phi_{n,i}^2). 
\] 

Now, consider the conditional expectation in the definition of $Z_{n,l}$. Notice that we have,
\begin{eqnarray*}
\sum_{i=j}^n \phi_{n,i}^2 &=& E_{\bar{\bbP}}[\sum_{i=j}^n \phi_{n,i}^2 | \sF_{n,j-1}]\\
&=& E_{\bar{\bbP}}[\sum_{i=j}^n \phi_{n,j}^2 | \sF_{n,j-1}]\\
&=& (n-j+1) E[\phi^2_{n,j} | \sF_{n,j-1}].
\end{eqnarray*}
The penultimate equation follows from the fact that  $\phi_{n,l}^2$ are exchangeable, from which the last equality follows by the fact that they are also identically distributed. This implies that
\[
E_{\bar \bbP}[\phi^2_{n,j} | \sF_{n,j-1}] = \frac{1}{n-j+1} \sum_{i=j}^n \phi_{n,i}^2
\]
Thus, the martingale difference sequence $Z_{n,l}$ has the compact representation 
\[
Z_{n,l} = \phi_{n,l}^2 - \frac{1}{n-l+1} \sum_{i=l}^n \phi_{n,i}^2.
 \]

In order to obtain a bound on $\sum_{l=1}^n \phi_{n,l}^2 - 1$, we first bound $\sum_{l=1}^n \frac{1}{n-l+1} \sum_{i=l}^n \phi_{n,i}^2$ from above. Recall that $\phi_{n,l} \leq 2T/\sqrt{n}$, so that

\begin{eqnarray*}
\sum_{l=1}^n \frac{1}{n-l+1} \sum_{i=l}^n \phi_{n,i}^2  &=& \sum_{l=0}^{n-1} \frac{1}{n-l} \sum_{j=l+1}^n \phi^2_{n,j}\\
&\leq& \frac{4T^2}{n} \sum_{l=0}^{n-1} \left( 1 - \frac{1}{n-l}
\right) \\
&\approx& 4 T^2 - 4T^2 \frac{\log n}{n} + \frac{o(n)}{n}.
\end{eqnarray*}

Therefore, it follows that 
\[
\sum_{l=1}^n Z_{n,l} \geq \sum_{l=1}^n \phi_{n,l}^2 - 4T^2 + 4T^2 \frac{\log n}{n} - \frac{o(n)}{n}.
\]
Now, consider the event $\{\omega \in \bar{\Omega} | \sum_{l=1}^n \phi_{n,l}^2 -1 \geq \e \}$, where $\e > 0$. From the bound, it follows that
\[
\sum_{l=1}^n Z_{n,l} \geq \sum_{l=1}^n \phi_{n,l}^2 - 4T^2 + 4T^2 \frac{\log n}{n}  \geq \e +1 - 4T^2 + 4T^2 \frac{\log n}{n}.
\]

Using the Azuma-Hoeffding inequality as described above, it follows that
\[
\bar{\bbP}(\sum_{l=1}^n Z_{n,l} \geq \e +1 - 4T^2 + 4T^2 \frac{\log n}{n}) \leq \exp \left( - \frac{(\e +1 - 4T^2 + 4T^2 \frac{\log n}{n})^2}{n \times \frac{64 T^4}{n^2}} \right).
\]
(ignoring the $o(n)$ term). The bound in the numerator on the RHS follows by the fact that
\[
|Z_{n,l}| \leq |\phi_{n,l}^2| + \frac{1}{n-l+1} \sum_{j=l}^n |\phi_{n,j}^2| \leq 2 |\phi_{n,l}^2|,
\]
and noting that $\phi_{n_l}^2 \leq 2 T^2/n$. Considering the expression being exponentiated on the RHS, it is a straightforward exercise to see that as $n \to \infty$, the expression tends to $\infty$, in turn implying that 
\[
\bar{\bbP}(\sum_{l=1}^n \phi_{n,l}^2 - 1 \geq \e) \leq \bar{\bbP}(\sum_{l=1}^n Z_{n,l} \geq \e +1 - 4T^2 + 4T^2 \frac{\log n}{n}) \to 0,
\]
as $n \to \infty$.

Next, note that since $\phi_{n,l}^2 \geq 0$ for all $1 \leq l \leq n$, it follows that $\sum_{l=1}^n Z_{n,l} \leq \sum_{n=1}^l \phi_{n,l}^2$. Clearly,
\[
\bar{\bbP}(\sum_{l=1}^n \phi_{n,l}^2 < 1 - \e) \leq \bar{\bbP}(\sum_{l=1}^n Z_{n,l} < 1 - \e). 
\] 
Using the Azuma-Hoeffding inequality again, we have
\[
\bar{\bbP}(\sum_{l=1}^n Z_{n,l} < 1 - \e) \leq \exp\left(-\frac{(1-\e)^2}{n \times \frac{64 T^4}{n^2}} \right).
\]
As $n \to \infty$, the LHS tends to 0, exponentially fast. Thus, it
follows that $|\sum_{l=1}^n \phi_{n,l}^2 - 1| \stackrel{P}{\to} 0$ as
$n \to \infty$, completing the proof.
\EndProof

\subsection{Proof of Theorem \ref{thm:counting}}
We first state a couple of useful lemmata. We will find it useful to work with a relaxed form
of $\bar N_n$: 

\begin{definition}[Relaxed Counting Process]
~$\tilde N_n(t) := \sup \{p \in [0,1] | \bar S_n(p) \leq t\}$.
\end{definition}
This process has the useful description as
the fraction of arrivals by time $t$. Clearly both $\bar N_n, \tilde
N_n \in \sD_{\lim}[0,1]$, and are asymptotically close as the following
proposition shows. $\| \cdot \|$ refers to the $\sup$ norm over
$[0,1]$. In the remainder of the section we work with the process $\tilde
N_n$. 

\begin{lemma} \label{lem:counting-counting}
~$\| \bar N_n - \tilde N_n\| \to 0$ as $n \to \infty$.
\end{lemma}
\Proof
Fix $p \in [0,1)$. Clearly, 
\[
\frac{m}{n} \leq p < \frac{m+1}{n},
\]
for some $m = 0, 1, \ldots, n-1$. Suppose that  $t \in
[0,1]$,
\[
|\bar N_n(t) - \tilde N_n(t)| = \tilde N_n(t) - \frac{m}{n},
\]
for some $m$ that is $t$ dependent. Using the upper bound we have
\(
|\bar N_n(t) - \tilde N_n(t)| < \frac{1}{n},
\) 
which is independent of $t$. The conclusion follows.
\EndProof

To complete the analysis we also require a definition of the
inverse function of $\bar S_n$.
\begin{definition}[Partial Sum Inverse]
\[
\bar S_n^{-1} (t) := \inf \{p \in [0,1] | \tilde N_n(p) > t\}.
\]
\end{definition}
The partial sum inverse and the relaxed counting process are related
by the expression: 
\(
\bar S_n^{-1}(t) = \tilde N_n( \tilde N_n^{-1} ( \tilde N_n(t))),
\)
where 
\(
\tilde N_n^{-1}(t) := \inf \{p \in [0,1] | \tilde N_n(p) > t\}.
\)
Clearly, the partial sum inverse and the counting process must be close,
asymptotically (as $\bar S_n$ converges to a continuous process in the
limit). The following lemma shows that this is indeed the case. 

\begin{lemma} \label{lem:counting-inverse}
\noindent (i)
\(
\|\bar N_n - \bar S_n^{-1}\| \to 0 
\)
as $n \to \infty$.\\
\noindent (ii)
\(
\sn (\bar N_n - \bar S_n^{-1}) \to 0
\)
as $n \to \infty$.
\end{lemma}
\Proof
\noindent (i) Fix $t \in [0,1]$. By definition it follows that $\bar S_n (\tilde
N_n(t)) \leq t$ and $\bar S_n (\bar S_n^{-1}(t)) > t$ (and $\bar S_n
(\bar S_n^{-1}(t)-) \leq t$). Thus, for any $\e > 0$, $\bar S_n (\tilde
N_n(t) + \e) > t$. In particular, $\tilde N_n(t) + \frac{1}{n} \geq
\bar S_n^{-1}(t)$. Since $\bar S_n$ is non-decreasing (since the
increments $\xi_{n,l} \geq 0$), it follows that
\[
\frac{1}{n} \geq \bar S_n^{-1}(t) - \tilde N_n(t) \geq 0,
\]
where the last inequality follows by definition. Combining this
expression with Lemma \ref{lem:counting-counting}, the conclusion
follows.

\noindent(ii) The result is an obvious corollary of the argument for
part $(i)$.
\EndProof

We now present the proof of Theorem \ref{thm:counting}.

By an application of Theorem 7.8.1 of
\cite{Wh01b}, the fSLLN in Theorem \ref{thm:conditioned-fslln} implies the convergence of the
corresponding inverse function, $\bar S_n^{-1}$ to $e$, and by part (i) of Lemma \ref{lem:counting-inverse} the
convergence of the counting process $\bar N_n$. This proves part (i)
of Theorem \ref{thm:counting}.

Next, note that
\[
\frac{1}{\sn} \sum_{l=1}^{\lfloor nt \rfloor } \xi_{n,l} - \sn t =
\frac{1}{\sn} \sum_{l=1}^{\lfloor nt \rfloor} \left( \xi_{n,l} - \mu_n
\right) + \sn \left( \lfloor nt \rfloor \mu_n - t \right) \Rightarrow W^0(t),
\]
u.o.c. of $[0,\infty)$ a.s. as $n \to \infty$, where the latter term
of the second expression converges to 0:
\[
\sn \left( \lfloor nt \rfloor \mu_n - t \right) \to 0
\]
as $n \to \infty$. Theorem 7.8.2 of \cite{Wh01b} implies that the fCLT
above implies the convergence of the scaled and centered inverse
process, and hence the counting process by part (ii) of Lemma \ref{lem:counting-inverse}
\EndProof

\subsection{Proof of Proposition \ref{prop:scheduled-avg-dist-limit}}
Fix $n \geq 1$  and $t \in [0,1]$. Define  $j_n^{*} := \inf \{j = 0, \ldots, n : t \in
(j/n-T, j/n + T]\}$ as the first arrival index such that $t$ is in the
support of $F_{n,j_n^{*}}$, and $j_n^{**} := \sup \{j = 0, \ldots, n : t
\geq j/n -T\}$ as the largest index such that $t$ is greater than the
lower bound of the support of $F_{n,j_n^{**}}$.

By the definition of the infimum and supremum, for any $\e > 0$,
\begin{eqnarray}
\frac{j_n^{*}}{n} - \e &< t- T  \leq& \frac{j_n^*}{n} ~~\text{and}\\
\frac{j_n^{**}}{n} &< t+T \leq& \frac{j_n^{**}}{n} + \e,
\end{eqnarray}
so that as $n \to \infty$, 
\(
\frac{j_n^*}{n} \to t-T
\)
~and~
\(
\frac{j_n^{**}}{n} \to t+T.
\)

Now, let $t \in [-T,T]$, then using the definition of $j_n^*$ and $j_n^{**}$,
\begin{eqnarray*}
\frac{1}{n} \sum_{i=1}^{n+1} F\left( t - \frac{i-1}{n}\right) &=&
\frac{1}{n} \sum_{i=1}^{j_n^{**}} F\left( t - \frac{i-1}{n}\right) \\
&=& \frac{1}{n} \sum_{i=1}^{j_n^{**}} \frac{t - \frac{i-1}{n} + T}{2T}
\\
&=& \frac{1}{2T} \left\{ (t+T) \frac{j_n^* - 1}{n} - \frac{1}{n^2}
  \frac{j_n^{**} (j_n^{**} - 1)}{2} \right\}\\
&\to& \frac{(t+T)^2}{4T},
\end{eqnarray*}
where the second equality follows by the fact that $j_n^* = 0$ since
$t \in [-T,T)$, and $j_n^{**} < n$ follows from the fact that $T \in
[0,0.5]$. The limit follows by the limit argument presented above for
$j_n^{**}$.

Next, fix $t \in (1-T,1+T]$. In this case, $j_n^{**} = n$ and we have
\begin{eqnarray*}
\frac{1}{n} \sum_{i=1}^{n+1} F\left( t - \frac{i-1}{n}\right) &=&
\frac{1}{n} \sum_{i=j_n^{*}}^{n+1} \frac{t+T - \frac{i-1}{n}}{2T} +
\frac{j_n^{*}}{n}\\
&=& \frac{1}{2T} \left \{ (t+T) \frac{n+1 - j_n^*}{n} - \frac{1}{n^2}
  \left[ \frac{(n+1)(n+2)}{2} - \frac{j_n^*(j_n^* + 1)}{2} \right]
\right \} + \frac{j_n^*}{n}\\
&\to& \frac{t+T}{2 T} - \frac{t^2 - T^2}{2 T} + \frac{(t - T)^2}{4 T} - \frac{1}{4 T} + (t-T),
\end{eqnarray*}
where the first equality follows by the fact that for all indices $j <
j_n^*$ the value of the distribution function at this $t$ is $1$ (in
effect, this $t$ is outside the support of these arrival
distributions). The limit, of course, is straightforward from those of
$j_n^*$.

Finally, for $t \in [T,1-T]$ we have 
\begin{eqnarray*}
\frac{1}{n} \sum_{i=1}^{n+1} F\left( t - \frac{i-1}{n}\right) &=&
\frac{j_n^*}{n} + \frac{1}{n} \sum_{i=j_n^*}^{j_n^{**}} \frac{t + T -
  \frac{i-1}{n}}{2 T}\\
&=& \frac{j_n^*}{n} + \left( \frac{t+T}{2T} \right) \left(
  \frac{j_n^{**} - j_n^*}{n} \right) - \frac{1}{2Tn^2} \left(
  \frac{j_n^{**}(j_n^** + 1)}{2} - \frac{j_n^*(j_n^*-1)}{2} \right)\\
&\to& t,
\end{eqnarray*}
where first equality follows from the fact that, for large enough $n$,
$j_n^* > 1$ and $j_n^** < n$. The rest of the argument is a
consequence of the limits for $j_n^*$ and $j_n^{**}$.

This completely characterizes the limit for a fixed $t$. Uniform
convergence can be shown by a similar argument as in the proof of Lemma \ref{lem:avg-dist-limit}.
\EndProof

\subsection{Proof of Corollary \ref{cor:scheduled-avg-dist-limit}}
Fix $t \in [-T,T)$. Notice that only users indexed by $p \in [0,t+T]$
can possibly arrive in this interval. Therefore, for such a $t$, we
have
\begin{eqnarray*}
\int_0^1 F_p(t) m(dp) &=& \int_0^{t+T} \left( \frac{t-p + T}{2T}
\right) dp \\
&=& \frac{(t+T)^2}{4 T}.
\end{eqnarray*}
This is matches the expression found in Proposition
\ref{prop:scheduled-avg-dist-limit}. Next, let $t \in [T,1-T]$. In
this instance, users indexed by $p \in [0,t-T)$ and $p \in (t+T,1]$ cannot
arrive in this interval. However, those indexed in the former interval
will have arrived by time $t$. This implies that
\begin{eqnarray*}
\int_0^1 F_p(t) m(dp) &=& \int_0^{t-T} m(dp) + \int_{t-T}^{t+T}
\frac{t- p + T}{2T} dp\\
&=& t,
\end{eqnarray*}
agreeing with Proposition
\ref{prop:scheduled-avg-dist-limit}. Finally, for $t \in (1-T,1+T]$,
only arrivals with index $p \in [t-T,1]$ can arrive in this interval
and, furthermore, all other arrivals will have happened by time
$t$. This implies that,
\begin{eqnarray*}
\int_0^1 F_p(t) m(dp) &=& \int_0^{t-T} dp + \int_{t-T}^1
\frac{t-p+T}{2T} dp\\
&=& \frac{t+T}{2 T} - \frac{t^2 - T^2}{2 T} + \frac{(t - T)^2}{4 T} - \frac{1}{4 T} + (t-T).
\end{eqnarray*}
This completes the proof.
\EndProof

\end{APPENDICES}

\section*{Acknowledgements}
The first author would  like to thank Vijay G. Subramanian for his help in proving a part of Theorem \ref{thm:condition-mean-rate} and for many insightful conversations. The authors would also like to thank Jim Dai, Peter Glynn, Bill Massey, Jamol Pender, Kavitha Ramanan, Sheldon Ross and Jean Walrand for helpful discussions and comments over the course of working through this paper.
\bibliographystyle{ormsv080}
\bibliography{../Journal/refs-queueing}
\end{document}